\documentclass[12pt,a4paper,twoside]{article}
\usepackage[T2A]{fontenc}
\usepackage[english,russian]{babel}

\fontsize{12}{21} \selectfont

\newtheorem{theorem}{Theorem}[section]

\newtheorem{definition}{Definition}[section]

\begin{document}

\noindent

\makeatletter
\renewcommand{\@evenhead}{\hfil
{\bf Andrei I. Bodrenko.}}
\renewcommand{\@oddhead}{\hfil
\small{\bf \underline{The solution of the Minkowski problem for
open surfaces in Riemannian space.}} }

\noindent
\bigskip
\begin{center}
{\large \bf The solution of the Minkowski problem for open
surfaces in Riemannian space.}
\end{center}
\medskip

\begin{center}
{\bf Andrei~I.~Bodrenko} \footnote{\copyright  Andrei~I.~Bodrenko,
associate professor,
Department of Mathematics,\\
Volgograd State University,
University Prospekt 100, Volgograd, 400062, RUSSIA.\\
E.-mail: bodrenko@mail.ru \qquad \qquad http://www.bodrenko.com}
\end{center}

\begin{center}
{\bf Abstract}\\
\end{center}

{\small Author reduces the Minkowski problem to the problem of
construction the G-deformations preserving the product of
principal curvatures for every point of surface in Riemannian
space. G-deformation transfers every normal vector of surface in
parallel along the path of the translation for each point of
surface. The continuous G-deformations preserving the product of
principal curvatures of surface with boundary are considered in
this article. The equations of deformations which are obtained in
this paper reduce to the nonlinear boundary-value problem. The
method of construction continuous G-deformations preserving the
product of principal curvatures  of surface with boundary and
 its qualitative analysis are presented in this article.}

\section*{Introduction}

\bigskip

The Minkowski problem (MP) is well known fundamental problem of
differential geometry. There have been published a number of
articles on this subject since 1903. But all authors have studied
this problem in Euclidean or pseudo-Euclidean spaces. They were H.
Minkowski, A.V. Pogorelov, A.D. Aleksandrov, W.J. Firey and many
others. We know many \\ generalizations of the MP in Euclidean and
pseudo-Euclidean spaces.

V.T. Fomenko [21] studied the MP by the methods of deformation
theory in Euclidean space.

The MP in Riemannian space differs substantially from the MP in
Euclidean space. The MP in Riemannian space is much more
complicate than the MP in Euclidean space.

Author of this article created for the first time the method of
finding solutions of the MP in Riemannian space using deformation
theory.

Author have studied AG-deformations in Euclidean spaces in [1-16].
It was very hard to obtain the equation system of AG-deformations
in Riemannian space and much harder to solve it. V.T. Fomenko in
[24] studied infinitesimal ARG-deformations in Riemannian space.
The methods developed in [1-16] are very useful for finding the
solution of the Minkowski problem.

\section*{\S 1. Basic definitions. Statement of the main result.}

Let $R^{3}$ be  the three-dimensional Riemannian space with metric
tensor $\tilde{a}_{\alpha\beta},$
  $F$ be the two-dimensional simply connected oriented
surface in $R^{3}$ with the boundary $\partial F.$

Let $F\in C^{m,\nu}, \nu \in (0;1) , m\ge 4,$ $\partial F\in
C^{m+1,\nu}.$ Let $F$ has all strictly positive principal
curvatures $k_{1}$ and $k_{2}$. Let F be oriented so that mean
curvature $H$ is strictly positive. Denote $K=k_{1}k_{2}.$

Let $F$ be given by immersion of the domain $D\subset E^{2}$ into
$R^{3}$ by the equation: $y^{\sigma}=f^{\sigma}(x), x\in D$, $f:D
\rightarrow R^{3}.$ Denote by $d\sigma(x)=\sqrt{g} dx^{1}\wedge
dx^{2}$ the area element of the surface $F$. We identify the
points of immersion of surface $F$ with the corresponding
coordinate sets in $R^{3}$. Without loss of generality we assume
that $D$ is unit disk. Let $x^{1}, x^{2}$ be the Cartesian
coordinates.

Symbol $_{,i}$ denotes covariant derivative in metric of surface
$F.$ Symbol $\partial_{i}$ denotes partial derivative by variable
$x^{i}.$ We will assume $\dot{f}\equiv \frac{d f}{d t}.$ We define
$\Delta(f)\equiv f(t)-f(0).$

Indices denoted by Greek alphabet letters define tensor
coordinates in Riemannian space $R^{3}.$ We use the following
rule: a formula is valid for all admissible values of indices if
there are no instructions for which values of indices it is valid.
We use the Einstein rule. We assume that integer $m_{1}$ satisfies
the condition $0 \leq m_{1}\leq m-2.$

We consider continuous deformation of the surface $F$: $\{F_{t}\}$
defined by the equations
$$y^{\sigma}_{t}=y^{\sigma} + z^{\sigma}(t), z^{\sigma}(0)\equiv 0,
 t\in[0;t_{0}], t_{0}>0. \eqno(1.1)$$

\begin{definition}
\label{definition 1}. Deformation $\{F_{t}\}$ is called the
continuous deformation preserving the product of principal
curvatures  ( or $M-$deformation [21]) if the following condition
holds: $\Delta(K)=0$ and $z^{\sigma}(t)$ is continuous by $t.$

\end{definition}

The deformation $\{F_{t}\}$ generates the following set of paths
in $R^{3}$
$$u^{\alpha_{0}}(\tau)=(y^{\alpha_{0}}+z^{\alpha_{0}}(\tau)), \eqno(1.2)$$
where $z^{\alpha_{0}}(0)\equiv 0, \tau \in [0;t], t\in[0;t_{0}],
t_{0}>0.$

\begin{definition}
\label{definition 2}. The deformation $\{F_{t}\}$ is called the
$G-$deformation if every normal vector of surface transfers in
parallel along the path of the translation for each point of
surface.
\end{definition}

Let, along the $\partial F$, be given vector field tangent to $F.$
We denote it by the following formula:
$$v^{\alpha}=l^{i}y^{\alpha}_{,i}.  \eqno(1.3)$$

We consider the boundary-value condition:
$$\tilde{a}_{\alpha\beta}z^{\alpha}v^{\beta}=\tilde{\gamma}(s,t) ,
s\in \partial D. \eqno(1.4)$$

Let $v^{\alpha}$ and $\tilde{\gamma}$ be of class $C^{m-2,\nu}.$

We denote:
$$\tilde{\lambda}_{k}=\tilde{a}_{\alpha\beta}y^{\alpha}_{,k}v^{\beta},
k=1,2.  \eqno(1.5)$$
$$\lambda_{k}=\frac{ \tilde{\lambda}_{k}}
{(\tilde{\lambda}_{1})^{2}+(\tilde{\lambda}_{2})^{2}}, k=1,2.
\eqno(1.6)$$

$$\lambda(s)=\lambda_{1}(s)+i \lambda_{2}(s), s\in \partial D. \eqno(1.7)$$

Let $n$ be the index of the given boundary-value condition
$$n=\frac{1}{2\pi}\Delta_{\partial D} \arg \lambda(s).  \eqno(1.8)$$

\begin{theorem}
\label{theorem1}. Let $F\in C^{m,\nu}, \nu \in (0;1) , m\ge 4,$
$\partial F \in C^{m+1,\nu}.$ Let $\tilde{a}_{\alpha\beta} \in
C^{m,\nu},$ $\exists M_{0}=const>0$ such that
$\|\tilde{a}_{\alpha\beta}\|_{m,\nu}<M_{0},$ $\|\partial
\tilde{a}_{\alpha\beta}\|_{m,\nu}<M_{0},$ $\|\partial^{2}
\tilde{a}_{\alpha\beta}\|_{m,\nu}<M_{0}.$ Let $v^{\beta},
\tilde{\gamma} \in C^{m-2,\nu}(\partial D),$ $\tilde{\gamma}$ is
continuously differentiable by $t.$ Let, at the point
$(x^{1}_{(0)},x^{2}_{(0)})$ of the domain $D,$ the following
condition holds: $\forall t : z^{\sigma}(t)\equiv 0.$

Then the following statements hold:

1) if $n > 0$ then there exist $t_{0}>0$ and
$\varepsilon(t_{0})>0$ such that for any admissible
$\tilde{\gamma}$ satisfying the condition:
$|\dot{\tilde{\gamma}}|_{m-2,\nu}\leq \varepsilon$
 for all $t\in [0, t_{0})$
there exists $(2n-1)-$parametric $MG-$deformation of class
$C^{m-2,\nu}(\bar{D})$ continuous by $t.$

2) if $n < 0$ then there exist $t_{0}>0$ and
$\varepsilon(t_{0})>0$ such that for any admissible
$\tilde{\gamma}$ satisfying the condition:
$|\dot{\tilde{\gamma}}|_{m-2,\nu}\leq \varepsilon(t_{0})$
 for all $t\in [0, t_{0})$
there exists nor more than one $MG-$deformation of class
$C^{m-2,\nu}(\bar{D})$ continuous by $t.$ 

3) if $n = 0$ then there exist $t_{0}>0$ and
$\varepsilon(t_{0})>0$ such that for any admissible
$\tilde{\gamma}$ satisfying the condition:
$|\dot{\tilde{\gamma}}|_{m-2,\nu}\leq \varepsilon$
 for all $t\in [0, t_{0})$
there exists one $MG-$deformation of class $C^{m-2,\nu}(\bar{D})$
continuous by $t.$
\end{theorem}

We denote:\\
$\tilde{a}_{\alpha\beta}(t)\equiv
\tilde{a}_{\alpha\beta}(y^{\sigma}+z^{\sigma}(t)),$
$\tilde{a}_{\alpha\beta}(0)\equiv \tilde{a}_{\alpha\beta}.$ \\
$\Gamma_{\alpha \beta}^{\gamma}(t)\equiv
 \Gamma_{\alpha \beta}^{\gamma}(y^{\sigma}+z^{\sigma}(t)),$
$\Gamma_{\alpha \beta}^{\gamma}(0)\equiv
 \Gamma_{\alpha \beta}^{\gamma}.$ \\
$b_{ij}(t)\equiv b_{ij}(y^{\sigma}+z^{\sigma}(t)),$
$b_{ij}(0)\equiv b_{ij}.$ $b(t)\equiv
b(y^{\sigma}+z^{\sigma}(t)),$
$b(0)\equiv b.$\\
$g_{ij}(t)\equiv g_{ij}(y^{\sigma}+z^{\sigma}(t)),$
$g_{ij}(0)\equiv g_{ij}.$ $g(t)\equiv
g(y^{\sigma}+z^{\sigma}(t)),$
$g(0)\equiv g.$ \\
$a^{j}(t)\equiv a^{j}, c(t)\equiv c,$ $z^{\sigma}(t)\equiv
z^{\sigma}.$

\section*{\S 2. Deduction of $G-$deformation formulas for surfaces in
Riemannian space.}

The deformation $\{F_{t}\}$ of surface $F$ is defined by (1.1).

We denote:
$$z^{\sigma}(t)=a^{j}(t)y^{\sigma},_{j} + c(t)n^{\sigma},  \eqno(2.1)$$
where $a^{j}(0)\equiv 0, c(0)\equiv 0.$ Therefore the deformation
of surface $F$ is defined by functions $a^{j}$ and $c.$

We denote:
$$
\nabla^{*}_{i} z^{\alpha}= z^{\alpha}_{,i}+\Gamma_{\gamma
\sigma}^{\alpha}(0) z^{\sigma}y^{\gamma},_{i}. \eqno(2.2)$$

Then we have (see [27]):
$$\nabla^{*}_{j}z^{\alpha}=(a^{i}b_{ij}+
c,_{j})n^{\alpha}+(a^{i},_{j}-cb_{jm}g^{mi})y^{\alpha},_{i}
\eqno(2.3)$$

The condition of $G-$deformation is equivalent to the following
equality: \\
$$\tilde{a}_{\alpha_{0} \beta}A^{*\alpha_{0}}_{i}(t) n^{\beta}=0, \eqno(2.4)$$
 where $A^{*\alpha_{0}}_{i}(t)$ is the result of parallel translation
of tensor $(y^{\alpha_{0}}_{,i}+z^{\alpha_{0}}_{,i}(t))$ from the
point
 $(y^{\alpha_{0}}+z^{\alpha_{0}}(t))$ to the point
  $(y^{\alpha_{0}})$ along the path of the translation for each
  point of surface by deformation, i.e. along the following curve:
$$u^{\alpha_{0}}(\tau)=(y^{\alpha_{0}}+z^{\alpha_{0}}(\tau)), \eqno(2.5 ) $$
where $z^{\alpha_{0}}(0)\equiv 0, \tau \in [0;t].$

Then parallel translation of tensor
 $y^{\alpha_{0}}_{,i}+z^{\alpha_{0}}_{,i}(t)$
 reduces to the following Cauchy problem:
 find $A^{\alpha_{0}}_{i}(t,\tau)$ satisfying the equations:
$$\frac{d A^{\alpha_{0}}_{i}(t,\tau)}{d\tau}+
\Gamma ^{\alpha_{0}}_{\beta\gamma}(\tau)
\dot{z}^{\beta}(\tau)A^{\gamma}_{i}(t,\tau)=0, \alpha_{0}=1,2,3,
\tau \in [0;t], \eqno(2.6)$$ with initial boundary value
conditions:
$$A^{\alpha_{0}}_{i}(t,t)= y^{\alpha_{0}}_{,i}+z^{\alpha_{0}}_{,i}(t).
\eqno(2.7)$$

We denote:
$$ A^{\alpha_{0}}_{(1)\gamma_{0}}(t,\tau)=
\int\limits_{\tau}^{t} \Gamma
^{\alpha_{0}}_{\beta_{0}\gamma_{0}}(\tau_{0})
\dot{z}^{\beta_{0}}(\tau_{0}) d\tau_{0}. \eqno(2.8) $$ For $k\ge
2$ we denote:
$$ A^{\alpha_{0}}_{(k)\gamma_{k-1}}(t,\tau)=
\int\limits_{\tau}^{t} \int\limits_{\tau_{0}}^{t} ...
\int\limits_{\tau_{k-2}}^{t}
 (\prod_{j=0}^{k-1} \Gamma ^{\alpha_{j}}_{\beta_{j}\gamma_{j}}(\tau_{j})
\dot{z}^{\beta_{j}}(\tau_{j}) )d\tau_{0}d\tau_{1}...d\tau_{k-1}.
\eqno(2.9)$$ and for $j\geq 1, k\geq 2 $ we assume
$$ \alpha_{j}\equiv \gamma_{j-1}. \eqno(2.10)$$

{\bf Lemma 2.1.} {\it The following inequalities hold:\\
$$\|A^{\alpha_{0}}_{(1)\gamma_{0}}(t,\tau) \|_{m_{1},\nu} \leq
t\|\Gamma\|^{(t)}_{m_{1},\nu} \|\dot{z}\|^{(t)}_{m_{1},\nu}, \tau
\in [0;t]. \eqno(2.11) $$
 $$\|A^{\alpha_{0}}_{(k)\gamma_{0}}(t,\tau) \|_{m_{1},\nu} \leq
t^{k}(\|\Gamma\|^{(t)}_{m_{1},\nu})^{k}
(\|\dot{z}\|^{(t)}_{m_{1},\nu})^{k}, \tau \in [0;t]. \eqno(2.12)
$$ } The proof follows from aspects of
$A^{\alpha_{0}}_{(k)\gamma_{0}}(t,\tau).$

{\bf Lemma 2.2.} {\it Let the following conditions hold:

1) metric tensor in $R^{3}$ satisfies the conditions: $\exists
M_{0}=const>0$ such that
$\|\tilde{a}_{\alpha\beta}\|_{m,\nu}<M_{0},$ $\|\partial
\tilde{a}_{\alpha\beta}\|_{m,\nu}<M_{0},$ $\|\partial^{2}
\tilde{a}_{\alpha\beta}\|_{m,\nu}<M_{0}.$

2) $\exists t_{0}>0$ such that $c(t), c_{,i}(t), a^{k}(t),
\partial_{i}a^{k}(t)$ are continuous by $t, \forall
t\in[0,t_{0}],$ $c(0)\equiv 0, c_{,i}(0)\equiv 0,
 a^{k}(0)\equiv 0, \partial_{i}a^{k}(0)\equiv 0.$

3) $\exists t_{0}>0$ such that $z^{\alpha}(t)\in C^{m-2,\nu} ,
z^{\alpha}_{,i}(t)\in C^{m-3,\nu},$ $\forall t\in[0,t_{0}].$

Then there exists $t_{*}>0$ such that $\forall t\in[0, t_{*}) $
exists the unique result $A^{*\alpha_{0}}_{i}(t)$
 of translation in parallel of
 the tensor $(y^{\alpha_{0}}_{,i}+z^{\alpha_{0}}_{,i}(t))$ from the point
 $(y^{\alpha_{0}}+z^{\alpha_{0}}(t))$ to the point
  $(y^{\alpha_{0}})$ along the path of the translation for each
  point of surface by deformation.
   $A^{*\alpha_{0}}_{i}(t)$ has the following representation
$$A^{*\alpha_{0}}_{i}(t)=y^{\alpha_{0}}_{,i}+z^{\alpha_{0}}_{,i}(t)+
\sum_{k=1}^{\infty}
(y^{\gamma_{k-1}}_{,i}+z^{\gamma_{k-1}}_{,i}(t))
A^{\alpha_{0}}_{(k)\gamma_{k-1}}(t,0) . \eqno(2.13) $$
$A^{*\alpha_{0}}_{i}(t)$ is of class $C^{m-3,\nu}$ and continuous
by $t.$ }

{\bf Proof.} Finding the result of translation in parallel of
 the tensor along the given curve
brings to the Cauchy problem of differential equation system.
Using the methods represented in [19, p. 56] we reduce the
differential equation system to the integral equation system which
is resolved by the method of successive approximations.

The null approximation is:
$$A^{*\alpha_{0}}_{(0)i}(t)=y^{\alpha_{0}}_{,i}+z^{\alpha_{0}}_{,i}(t).
\eqno(2.14) $$

The $p-$th $(p>0)$ approximation of Cauchy problem is:
$$A^{*\alpha_{0}}_{(p)i}(t)=y^{\alpha_{0}}_{,i}+z^{\alpha_{0}}_{,i}(t)+
\sum_{k=1}^{p} (y^{\gamma_{k-1}}_{,i}+z^{\gamma_{k-1}}_{,i}(t))
A^{\alpha_{0}}_{(k)\gamma_{k-1}}(t,0) . \eqno(2.15) $$ Taking into
account that $C^{m,\nu}$ is complete normed space,
 lemma 2.1. and using reasonings that are similar to the ones
from [19, p. 56] for solution of this Cauchy problem we get
 the proof of lemma 2.2.

{\bf Lemma 2.3.} {\it Let the conditions of lemma 2.2 hold. Then
there exists $t_{*}>0$ such that $\forall t\in[0, t_{*})$
the following holds:\\
$$\dot{A}^{*\alpha_{0}}_{i}(t)=\dot{z}^{\alpha_{0}}_{,i}(t)+
\sum_{k=1}^{\infty} \dot{z}^{\gamma_{k-1}}_{,i}(t)
A^{\alpha_{0}}_{(k)\gamma_{k-1}}(t,0)+ \sum_{k=1}^{\infty}
(y^{\gamma_{k-1}}_{,i}+z^{\gamma_{k-1}}_{,i}(t))
\dot{A}^{\alpha_{0}}_{(k)\gamma_{k-1}}(t,0). \eqno(2.16) $$
$\dot{A}^{*\alpha_{0}}_{i}(t)$ is of class $C^{m-3,\nu}$
 and continuous by $t.$
}

The proof follows from the rules of termwise differentiation of
functional series, lemmas 2.1., 2.2.  and the properties of space
$C^{m-3,\nu}.$

Let obtain the equations of $G-$deformation and transform them to
the appropriate for our method form.

We denote $\sigma \equiv \gamma_{k-1}.$ Then we have:
$$A^{*\alpha_{0}}_{i}(t)=y^{\alpha_{0}}_{,i}+z^{\alpha_{0}}_{,i}(t)+
(y^{\sigma}_{,i}+z^{\sigma}_{,i}(t))\sum_{k=1}^{\infty}
A^{\alpha_{0}}_{(k)\sigma}(t,0) . \eqno(2.17) $$

We denote: $$\nabla^{*}_{i}
z^{\alpha_{0}}(t)=z^{\alpha_{0}}_{,i}(t)+
\Gamma^{\alpha_{0}}_{\beta_{0} \sigma}(0) y^{\sigma}_{,i}
z^{\beta_{0}}(t). \eqno(2.18) $$

We can write (2.17) in the following form: \\
$$A^{*\alpha_{0}}_{i}(t)=y^{\alpha_{0}}_{,i}+z^{\alpha_{0}}_{,i}(t)+
\Gamma^{\alpha_{0}}_{\beta_{0} \sigma}(0) y^{\sigma}_{,i}
z^{\beta_{0}}(t)+$$
$$(y^{\sigma}_{,i}+z^{\sigma}_{,i}(t)) A^{\alpha_{0}}_{(1)\sigma}(t,0)-
\Gamma^{\alpha_{0}}_{\beta_{0} \sigma}(0) y^{\sigma}_{,i}
z^{\beta_{0}}(t)
+(y^{\sigma}_{,i}+z^{\sigma}_{,i}(t))\sum_{k=2}^{\infty}
A^{\alpha_{0}}_{(k)\sigma}(t,0) . \eqno(2.19) $$

Inserting (2.18) into (2.19), we get: \\
$$A^{*\alpha_{0}}_{i}(t)=y^{\alpha_{0}}_{,i}+\nabla^{*}_{i} z^{\alpha_{0}}(t)+
z^{\sigma}_{,i}(t) A^{\alpha_{0}}_{(1)\sigma}(t,0)+$$
$$y^{\sigma}_{,i} ( A^{\alpha_{0}}_{(1)\sigma}(t,0)-
\Gamma^{\alpha_{0}}_{\beta_{0} \sigma}(0) z^{\beta_{0}}(t))+
(y^{\sigma}_{,i}+z^{\sigma}_{,i}(t))\sum_{k=2}^{\infty}
A^{\alpha_{0}}_{(k)\sigma}(t,0) . \eqno(2.20) $$

Therefore we can write (2.20) as: \\
$$A^{*\alpha_{0}}_{i}(t)=y^{\alpha_{0}}_{,i}+\nabla^{*}_{i} z^{\alpha_{0}}(t)+
z^{\sigma}_{,i}(t)\sum_{k=1}^{\infty}
A^{\alpha_{0}}_{(k)\sigma}(t,0)+ $$
$$
y^{\sigma}_{,i} ( A^{\alpha_{0}}_{(1)\sigma}(t,0)-
\Gamma^{\alpha_{0}}_{\beta_{0} \sigma}(0) z^{\beta_{0}}(t))
+y^{\sigma}_{,i}\sum_{k=2}^{\infty}
A^{\alpha_{0}}_{(k)\sigma}(t,0) . \eqno(2.21) $$

We denote: \\
$$S^{\alpha_{0}}_{(1)\sigma}(t,0)=
 \sum_{k=1}^{\infty} A^{\alpha_{0}}_{(k)\sigma}(t,0). \eqno(2.22) $$
$$S^{\alpha_{0}}_{(2)\sigma}(t,0)=
 \sum_{k=2}^{\infty} A^{\alpha_{0}}_{(k)\sigma}(t,0). \eqno(2.23) $$

Using (2.22), (2.23) from (2.21), we obtain: \\
$$A^{*\alpha_{0}}_{i}(t)=y^{\alpha_{0}}_{,i}+\nabla^{*}_{i} z^{\alpha_{0}}(t)+
z^{\sigma}_{,i}(t) S^{\alpha_{0}}_{(1)\sigma}(t,0)+$$
$$
y^{\sigma}_{,i} ( A^{\alpha_{0}}_{(1)\sigma}(t,0)-
\Gamma^{\alpha_{0}}_{\beta_{0} \sigma}(0) z^{\beta_{0}}(t))
+y^{\sigma}_{,i} S^{\alpha_{0}}_{(2)\sigma}(t,0) . \eqno(2.24) $$

We denote: \\
$$S^{\alpha_{0}}_{(3) i}(t,0)=
y^{\sigma}_{,i} ( A^{\alpha_{0}}_{(1)\sigma}(t,0)-
\Gamma^{\alpha_{0}}_{\beta_{0} \sigma}(0) z^{\beta_{0}}(t))+
y^{\sigma}_{,i} S^{\alpha_{0}}_{(2)\sigma}(t,0) . \eqno(2.25) $$

From (2.24) and (2.25) we get: \\
$$A^{*\alpha_{0}}_{i}(t)=y^{\alpha_{0}}_{,i}+\nabla^{*}_{i} z^{\alpha_{0}}(t)+
z^{\sigma}_{,i}(t) S^{\alpha_{0}}_{(1)\sigma}(t,0)+
S^{\alpha_{0}}_{(3) i}(t,0) . \eqno(2.26) $$

From (2.1), we get:
$$z^{\sigma},_{i}(t)=a^{j},_{i}y^{\sigma},_{j}+c,_{i}n^{\sigma}+
a^{j}y^{\sigma}_{,j,i}+cn^{\sigma}_{,i}. \eqno(2.27)$$

Inserting (2.27) into (2.26), we have: \\
$$A^{*\alpha_{0}}_{i}(t)=y^{\alpha_{0}}_{,i}+\nabla^{*}_{i} z^{\alpha_{0}}(t)+
S^{\alpha_{0}}_{(3) i}(t,0)+
(a^{j},_{i}y^{\sigma},_{j}+c,_{i}n^{\sigma})
S^{\alpha_{0}}_{(1)\sigma}(t,0)+$$
$$
(a^{j}y^{\sigma}_{,j,i}+cn^{\sigma}_{,i})S^{\alpha_{0}}_{(1)\sigma}(t,0).
\eqno(2.28) $$

We denote: \\
$$S^{\alpha_{0}}_{(4) i}(t,0)=S^{\alpha_{0}}_{(3) i}(t,0)+
(a^{j}y^{\sigma}_{,j,i}+cn^{\sigma}_{,i})S^{\alpha_{0}}_{(1)\sigma}(t,0).
\eqno(2.29) $$

Insertion (2.29) into (2.28), we obtain: \\
$$A^{*\alpha_{0}}_{i}(t)=y^{\alpha_{0}}_{,i}+\nabla^{*}_{i} z^{\alpha_{0}}(t)+
S^{\alpha_{0}}_{(4) i}(t,0)+
(a^{j},_{i}y^{\sigma},_{j}+c,_{i}n^{\sigma})
S^{\alpha_{0}}_{(1)\sigma}(t,0). \eqno(2.30) $$

Consider the following formula: \\
$$a^{j},_{i}=\partial_{i}(a^{j})+\Gamma^{j}_{pi} a^{p}. \eqno(2.31) $$

Then, form (2.30) and (2.31), we have: \\
$$A^{*\alpha_{0}}_{i}(t)=y^{\alpha_{0}}_{,i}+\nabla^{*}_{i} z^{\alpha_{0}}(t)+
(\partial_{i}(a^{j})y^{\sigma},_{j}+c,_{i}n^{\sigma})
S^{\alpha_{0}}_{(1)\sigma}(t,0)+$$
$$ S^{\alpha_{0}}_{(4) i}(t,0)+
\Gamma^{j}_{pi}
a^{p}y^{\sigma},_{j}S^{\alpha_{0}}_{(1)\sigma}(t,0). \eqno(2.32)
$$

We denote: \\
$$S^{\alpha_{0}}_{(5) i}(t,0)=
S^{\alpha_{0}}_{(4) i}(t,0)+ \Gamma^{j}_{pi} a^{p}y^{\sigma},_{j}
S^{\alpha_{0}}_{(1)\sigma}(t,0). \eqno(2.33) $$
$$T^{\alpha_{0}}_{0}(t,0)=n^{\sigma} S^{\alpha_{0}}_{(1)\sigma}(t,0).
\eqno(2.34) $$
$$T^{\alpha_{0}}_{j}(t,0)=y^{\sigma},_{j} S^{\alpha_{0}}_{(1)\sigma}(t,0),
j=1,2. \eqno(2.35) $$

Using (2.33), (2.34) and (2.35), we  can write (2.32) in the
following
form: \\
$$A^{*\alpha_{0}}_{i}(t)=y^{\alpha_{0}}_{,i}+\nabla^{*}_{i} z^{\alpha_{0}}(t)+
c,_{i}T^{\alpha_{0}}_{0}(t,0)+
\partial_{i}(a^{j})T^{\alpha_{0}}_{j}(t,0)+S^{\alpha_{0}}_{(5) i}(t,0).
\eqno(2.36) $$

For $G-$deformation the following condition holds:\\
$$\tilde{a}_{\alpha_{0}\beta_{0}}A^{*\alpha_{0}}_{i}(t)n^{\beta_{0}}=0.
\eqno(2.37) $$

Insertion (2.36) into (2.37), we have: \\
$$\tilde{a}_{\alpha_{0}\beta_{0}}A^{*\alpha_{0}}_{i}(t)n^{\beta_{0}}=
a^{l}b_{li}+c,_{i}+ c,_{i} \tilde{a}_{\alpha_{0}\beta_{0}}
T^{\alpha_{0}}_{0}(t,0)n^{\beta_{0}}+$$
$$\partial_{i}(a^{j}) \tilde{a}_{\alpha_{0}\beta_{0}}
T^{\alpha_{0}}_{j}(t,0)n^{\beta_{0}}+
\tilde{a}_{\alpha_{0}\beta_{0}} S^{\alpha_{0}}_{(5)
i}(t,0)n^{\beta_{0}}. \eqno(2.38) $$

We denote: \\
$$N_{j}(t,0)=
\tilde{a}_{\alpha_{0}\beta_{0}}
T^{\alpha_{0}}_{j}(t,0)n^{\beta_{0}}, j=0,1,2,$$
$$Q_{i}(t,0)=
\tilde{a}_{\alpha_{0}\beta_{0}} S^{\alpha_{0}}_{(5)
i}(t,0)n^{\beta_{0}}. \eqno(2.39) $$

Then the equations of $G-$deformation are:\\
$$ a^{l}b_{li}+(1+N_{0}(t,0))c,_{i}+
\partial_{i}a^{j} N_{j}(t,0) +Q_{i}(t,0)=0,i=1,2. \eqno(2.40) $$

\section*{\S 3. The estimations of norms.}

We denote:

$\|S_{(p)}\|^{(t)}_{m_{1},\nu}=\max_{\alpha_{0},\sigma}
\|S^{\alpha_{0}}_{(p)\sigma}(t,0)\|_{m_{1},\nu}, p=1,2.$

$\|S_{(l)}\|^{(t)}_{m_{1},\nu}=
\max_{\alpha_{0},i}\|S^{\alpha_{0}}_{(l) i}(t,0)\|_{m_{1},\nu},
l=3,4,5.$

$\|T\|^{(t)}_{m_{1},\nu}=
\max_{\alpha_{0},j=0,1,2}\|T^{\alpha_{0}}_{j}(t,0)\|_{m_{1},\nu}.$

$\|N\|^{(t)}_{m_{1},\nu}=
\max_{j=0,1,2}\|N_{j}(t,0)\|_{m_{1},\nu}.$

$\|Q\|^{(t)}_{m_{1},\nu}= \max_{i}\|Q_{i}(t,0)\|_{m_{1},\nu}.$

{\bf Lemma 3.1.} {\it The following estimations hold:

1) $\|S_{(1)}\|^{(t)}_{m_{1},\nu}\leq \sum_{k=1}^{\infty}
(t\|\Gamma\|^{(t)}_{m_{1},\nu}\|\dot{z}\|^{(t)}_{m_{1},\nu})^{k}.$

2) $\|S_{(2)}\|^{(t)}_{m_{1},\nu}\leq \sum_{k=2}^{\infty}
(t\|\Gamma\|^{(t)}_{m_{1},\nu}\|\dot{z}\|^{(t)}_{m_{1},\nu})^{k}.$

3) $\|S_{(3)}\|^{(t)}_{m_{1},\nu}\leq
Kt\|z\|^{(t)}_{m_{1},\nu}\|\dot{z}\|^{(t)}_{m_{1},\nu}+
K_{2}\sum_{k=2}^{\infty}
(t\|\Gamma\|^{(t)}_{m_{1},\nu}\|\dot{z}\|^{(t)}_{m_{1},\nu})^{k}.
$

4) $\|S_{(4)}\|^{(t)}_{m_{1},\nu}\leq
\|S_{(3)}\|^{(t)}_{m_{1},\nu}+
M_{9}\|z\|^{(t)}_{m_{1},\nu}\|S_{(1)}\|^{(t)}_{m_{1},\nu}.$

5) $\|S_{(5)}\|^{(t)}_{m_{1},\nu}\leq
\|S_{(4)}\|^{(t)}_{m_{1},\nu}+M_{10}\|S_{(1)}\|^{(t)}_{m_{1},\nu}.$

6) $\|T\|^{(t)}_{m_{1},\nu}\leq
M_{11}\|S_{(1)}\|^{(t)}_{m_{1},\nu}.$

7) $\|N\|^{(t)}_{m_{1},\nu}\leq
M_{14}\|S_{(1)}\|^{(t)}_{m_{1},\nu}.$

8) $\|Q\|^{(t)}_{m_{1},\nu}\leq
M_{15}\|S_{(5)}\|^{(t)}_{m_{1},\nu}.$ }

The proof of lemma follows from the forms of estimated functions
 and properties of norms in the space  $C^{m_{1},\nu}.$

{\bf Lemma 3.2.} {\it Let the following conditions hold:

1) metric tensor in $R^{3}$ satisfies the conditions: $\exists
M_{0}=const>0$ such that
$\|\tilde{a}_{\alpha\beta}\|_{m,\nu}<M_{0},$ $\|\partial
\tilde{a}_{\alpha\beta}\|_{m,\nu}<M_{0},$ $\|\partial^{2}
\tilde{a}_{\alpha\beta}\|_{m,\nu}<M_{0}.$

2) $\exists t_{0}>0$ such that $c(t), c_{,i}(t), a^{k}(t),
\partial_{i}a^{k}(t)$ are continuous by $t, \forall
t\in[0,t_{0}],$ $c(0)\equiv 0, c_{,i}(0)\equiv 0,
 a^{k}(0)\equiv 0, \partial_{i}a^{k}(0)\equiv 0.$

3) $\exists t_{0}>0$ such that $z^{\alpha}(t)\in C^{m-2,\nu} ,
z^{\alpha}_{,i}(t)\in C^{m-3,\nu},$ $\forall t\in[0,t_{0}].$

Then $\forall \varepsilon>0 \exists t_{0}>0$ such that

1) $ \|S_{(i)}\|^{(t)}_{m-2,\nu} \leq \varepsilon, \forall t\in
[0,t_{0}],  i=\overline{1,5}.$

2) $ \|T\|^{(t)}_{m-2,\nu} \leq \varepsilon, \forall t\in
[0,t_{0}].$

3) $ \|N\|^{(t)}_{m-2,\nu} \leq \varepsilon, \forall t\in
[0,t_{0}].$

4) $ \|Q\|^{(t)}_{m-2,\nu} \leq \varepsilon, \forall t\in
[0,t_{0}].$ }

The proof of lemma follows from the forms of considered functions
 and properties of the space $C^{m,\nu}$ and previous lemmas.

\section*{\S 4. Transformation of the $G-$deformations equations.}

We introduce conjugate isothermal coordinate system where
$b_{ii}=V,i=1,2, b_{12}=b_{21}=0.$ Then we have the equation
system from (2.40):
$$ c,_{1}(1+N_{0})+Va^{1}+N_{k}\partial_{1}a^{k}+Q_{1}=0$$
$$ c,_{2}(1+N_{0})+Va^{2}+N_{k}\partial_{2}a^{k}+Q_{2}=0 \eqno(4.1)$$

We differentiate the first equation by $x^{2},$ the second one by
$x^{1},$ and subtract from the first equation the second one. Then
we obtain:
$$ V\partial_{2}a^{1}-V\partial_{1}a^{2}+
c_{,1}\partial_{2}N_{0}-c_{,2}\partial_{1}N_{0}+
\partial_{1}a^{k}\partial_{2}N_{k}-
\partial_{2}a^{k}\partial_{1}N_{k}+ $$
$$ \partial_{2}Va^{1}-\partial_{1}Va^{2}+
\partial_{2}Q_{1}-\partial_{1}Q_{2}=0. \eqno(4.2)$$
We denote
$$
\Psi_{1}= -(c_{,1}\partial_{2}N_{0}-c_{,2}\partial_{1}N_{0}+
\partial_{1}a^{k}\partial_{2}N_{k}-
\partial_{2}a^{k}\partial_{1}N_{k}+
\partial_{2}Q_{1}-\partial_{1}Q_{2})/V. \eqno(4.3)
$$
Then, from (4.2) and (4.3), we have the following equation:
$$ \partial_{2}a^{1}-\partial_{1}a^{2}+p_{k}a^{k}=\Psi_{1}, \eqno(4.4) $$
where $p_{1}=\partial_{2}(\ln V), p_{2}=-\partial_{1}(\ln V).$
Note that $p_{k}$ do not depend on $t.$

Differentiating the equation (4.4) by $t$ we obtain the following
equation:
$$ \partial_{2}\dot{a}^{1}-\partial_{1}\dot{a}^{2}+
p_{k}\dot{a}^{k}=\dot{\Psi}_{1}. \eqno(4.5) $$

\section*{\S 5. Solution of the equation system (4.1):\\
 finding function $\dot{c}$ on functions $\dot{a}^{i}.$}

We will solve the equation system (4.1) assuming that functions
$a^{1}$ and $a^{2}$ are given. Note that $N_{k}, Q_{i}$ depend
only on $c , \dot{c} , a^{i} , \dot{a}^{i} .$ Function $V$ does
not depend on $c ,a^{i} $. We will use the following formulas.

$$c(x^{1},x^{2},t)= \int\limits_{0}^{t} \dot{c}(x^{1},x^{2},\tau)d \tau,
(c(x^{1},x^{2},0)=0). \eqno(5.1)$$

$$a^{i}(x^{1},x^{2},t)=
\int\limits_{0}^{t} \dot{a}^{i}(x^{1},x^{2},\tau)d \tau,
(a^{i}(x^{1},x^{2},0)=0). \eqno(5.2) $$

For functions $a^{i}_{,j}$ we will use the following formula:

$$a^{i}_{,k}(x^{1},x^{2},t)=
\int\limits_{0}^{t} \dot{a}^{i}_{,k}(x^{1},x^{2},\tau)d \tau,
(a^{i}_{,k}(x^{1},x^{2},0)=0). \eqno(5.3) $$

Formulas (5.1), (5.2) and (5.3) establish the connections between
functions $c , a^{i}$ and $\dot{c} , \dot{a}^{i} .$ It means that
if the functions $\dot{c} , \dot{a}^{i} $ are found then the
functions $c , a^{i} $ are found also.

Therefore we pass on to the new equation system (5.4) where there
we will consider functions $\dot{c} , \dot{a}^{i} ,\dot{c}_{,i} ,
\dot{a}^{i}_{,j} . $ We differentiate the equation system (4.1) by
$t$ and get (5.4). Note that $N_{k}, Q_{i},\dot{N}_{k},
\dot{Q}_{i}, $ depend only on $c , \dot{c} , a^{i} , \dot{a}^{i} $
and therefore depend only on $\dot{c} , \dot{a}^{i} .$ We can show
this by differentiating $N_{k}, Q_{i}$ by $t.$

Then we obtain equation system  for $\dot{c}$.
$$ \dot{c},_{1}=-\frac{d}{dt} \Biggl(
\frac{Va^{1}+N_{k}\partial_{1}a^{k}+Q_{1}}{(1+N_{0})} \Biggr)$$
$$ \dot{c},_{2}=- \frac{d}{dt} \Biggl(
\frac{Va^{2}+N_{k}\partial_{2}a^{k}+Q_{2}}{(1+N_{0})} \Biggr)
\eqno(5.4)$$

We can present equation system (5.4) as following:
$$ \dot{c},_{i}=-V\dot{a}^{i}
-\Biggl( \frac{-V\dot{a}^{i}N_{0}+N_{k}\partial_{i}\dot{a}^{k}+
\dot{N}_{k}\partial_{i}a^{k}+\dot{Q}_{i}}{(1+N_{0})}\Biggr)
 + \frac{ \dot{N_{0}}(Va^{i}+N_{k}\partial_{i}a^{k}+Q_{i})}{(1+N_{0})^{2}}
 \eqno(5.5)$$

 Then we transform (5.5) into integral equation relative to function
  $\dot{c}.$
Let $l^{*}$ be arbitrary admissible curve in $D$ starting at the
point
 $(x^{1}_{(0)},x^{2}_{(0)})$ and given by the equations
$x^{1}=x^{1}(s), x^{2}=x^{2}(s).$ Then we have the following
equation.

$$ \dot{c}(x^{1},x^{2},t)=$$
$$\int\limits_{(x^{1}_{(0)},x^{2}_{(0)})}^{(x^{1},x^{2})}
\Biggl(- \frac{-V\dot{a}^{1}N_{0}+N_{k}\partial_{1}\dot{a}^{k}+
\dot{N}_{k}\partial_{1}a^{k}+\dot{Q}_{1}}{(1+N_{0})}+ \frac{
\dot{N_{0}}(Va^{1}+N_{k}\partial_{1}a^{k}+Q_{1})}{(1+N_{0})^{2}}
\Biggr) d\tilde{x}^{1}+$$
$$
\Biggl(-\frac{-V\dot{a}^{2}N_{0}+N_{k}\partial_{2}\dot{a}^{k}+
\dot{N}_{k}\partial_{2}a^{k}+\dot{Q}_{2}}{(1+N_{0})}+ \frac{
\dot{N_{0}}(Va^{2}+N_{k}\partial_{2}a^{k}+Q_{2})}{(1+N_{0})^{2}}
\Biggr)d\tilde{x}^{2}+$$
$$\int\limits_{(x^{1}_{(0)},x^{2}_{(0)})}^{(x^{1},x^{2})}
\Biggl(-V\dot{a}^{1}\Biggr) d\tilde{x}^{1}+
\Biggl(-V\dot{a}^{2}\Biggr)d\tilde{x}^{2} \eqno(5.6)$$

Then the equation (5.6) along $l^{*}$ takes the form:

$$ \dot{c}(x^{1},x^{2},t)=$$
$$\int\limits_{0}^{s} \Biggl(
\Biggl(- \frac{-V\dot{a}^{1}N_{0}+N_{k}\partial_{1}\dot{a}^{k}+
\dot{N}_{k}\partial_{1}a^{k}+\dot{Q}_{1}}{(1+N_{0})} +\frac{
\dot{N_{0}}(Va^{1}+N_{k}\partial_{1}a^{k}+Q_{1})}{(1+N_{0})^{2}}
\Biggr){x^{1}}'(s1)+$$
$$\Biggl(-\frac{-V\dot{a}^{2}N_{0}+N_{k}\partial_{2}\dot{a}^{k}+
\dot{N}_{k}\partial_{2}a^{k}+\dot{Q}_{2}}{(1+N_{0})}
+\frac{\dot{N_{0}}(Va^{2}+N_{k}\partial_{2}a^{k}+Q_{2})}{(1+N_{0})^{2}}
\Biggr){x^{2}}'(s1) \Biggr) ds1+$$
$$\int\limits_{0}^{s}
\Biggl(- V\dot{a}^{1}(s1) {x^{1}}'(s1)
-V\dot{a}^{2}(s_{1}){x^{2}}'(s1)\Biggr) ds1. \eqno(5.7)$$

The equation (5.7) is nonlinear integral equation. We will show
that (5.7) has unique solution of class of continuous functions
for any continuous functions $\dot{a}^{i}$ and
$\partial_{p}\dot{a}^{i}.$

The equation (5.7) takes the form
$\dot{c}=L_{a}(\dot{c})+\gamma_{t},$ where operator $L_{a}$ has
explicit form.
$$\gamma_{t}(s)=\int\limits_{0}^{s}
\Biggl(- V\dot{a}^{1}(s1) {x^{1}}'(s1)
-V\dot{a}^{2}(s_{1}){x^{2}}'(s1)\Biggr)ds1. \eqno(5.8) $$

Therefore every pair of functions $\dot{a}^{i}\in C^{m-2,\nu}$
corresponds to the unique function $\dot{c}\in C^{m-2,\nu}$ and
therefore
to the unique function $c \in C^{m-2,\nu}:$\\
$c(x^{1},x^{2},t)= \int\limits_{0}^{t} \dot{c}(x^{1},x^{2},\tau)d
\tau, $ ($c(x^{1},x^{2},0)=0$).

Then the equation along $l^{*}$ takes the form:\\
$$ \dot{c}(x^{1},x^{2},t)=
\int\limits_{0}^{s} \Biggl(K_{a}(s_{1},\dot{c}(s_{1}))\Biggr)
ds_{1}+ \gamma_{t}(s), \eqno(5.9) $$ where
$$K_{a}(s_{1},\dot{c}(s_{1}))=$$
$$\Biggl(
\Biggl(- \frac{-V\dot{a}^{1}N_{0}+N_{k}\partial_{1}\dot{a}^{k}+
\dot{N}_{k}\partial_{1}a^{k}+\dot{Q}_{1}}{(1+N_{0})} +\frac{
\dot{N_{0}}(Va^{1}+N_{k}\partial_{1}a^{k}+Q_{1})}{(1+N_{0})^{2}}
\Biggr){x^{1}}'(s1)+$$
$$
\Biggl(-\frac{-V\dot{a}^{2}N_{0}+N_{k}\partial_{2}\dot{a}^{k}+
\dot{N}_{k}\partial_{2}a^{k}+\dot{Q}_{2}}{(1+N_{0})}
+\frac{\dot{N_{0}}(Va^{2}+N_{k}\partial_{2}a^{k}+Q_{2})}{(1+N_{0})^{2}}
\Biggr){x^{2}}'(s1) \Biggr). \eqno(5.10)$$

We denote: $$L_{a}(\dot{c})= \int\limits_{0}^{s}
\Biggl(K_{a}(s_{1},\dot{c}(s_{1}))\Biggr) ds_{1}.
 \eqno(5.11)$$

We will investigate the decidability problem of the equation in
the space  $C^{m-2,\nu}(\bar{D})$:
$$\dot{c}=L_{a}(\dot{c})+\gamma_{t}. \eqno(5.12) $$
We will solve equation (5.12) by the method of successive
approximations.

{\bf Lemma 5.1.} {\it Let the following conditions hold:

1) metric tensor in $R^{3}$ satisfies the conditions: $\exists
M_{0}=const>0$ such that
$\|\tilde{a}_{\alpha\beta}\|_{m,\nu}<M_{0},$ $\|\partial
\tilde{a}_{\alpha\beta}\|_{m,\nu}<M_{0},$ $\|\partial^{2}
\tilde{a}_{\alpha\beta}\|_{m,\nu}<M_{0}.$

2) $\exists t_{0}>0$ such that $a^{k}(t), \partial_{i}a^{k}(t),
\dot{a}^{k}(t), \partial_{i}\dot{a}^{k}(t)$ are continuous by $t,
\forall t\in[0,t_{0}],$ $ a^{k}(0)\equiv 0,
\partial_{i}a^{k}(0)\equiv 0.$

3) $\exists t_{0}>0$ such that $a^{i}(t)\in C^{m-2,\nu} ,
\partial_{k} a^{i}(t)\in C^{m-3,\nu},$ $\forall t\in[0,t_{0}].$

Then $\exists t_{*}>0$ such that the equation
 $\dot{c}=L_{a}(\dot{c})+\gamma_{t}$ $\forall t\in[0,t_{*}].$
 has unique solution of class $C^{m-2,\nu}$ continuous by $t.$
}

{\bf Proof.}

We construct the sequence of functions $\{\dot{c}^{(k)}\}$: we
find function $\dot{c}^{(1)}$ from the equation

$$ \dot{c}^{(1)}(x^{1},x^{2},t)=\int\limits_{0}^{s}
\Biggl(- V\dot{a}^{1}(s1) {x^{1}}'(s1)
-V\dot{a}^{2}(s_{1}){x^{2}}'(s1)\Biggr)ds1, \eqno(5.13)$$

we find function $\dot{c}^{(k)}, k>1 $ from the equation
$$\dot{c}^{(k)}=L_{a}(\dot{c}^{(k-1)})+\gamma_{t}. \eqno(5.14) $$
The sequence of functions $\{\dot{c}^{(k)}\}$ is determined
uniquely and functions $\dot{c}^{(k)}$ are of class
$C^{m-2,\nu}(\bar{D})$.

We will show that the sequence of functions $\{\dot{c}^{(k)}\}$ is
bounded in the space $C^{m-2,\nu}(\bar{D})$.

For any $\varepsilon>0$ there exists $t_{0}>0$ such that
 for all $t\in[0,t_{0})$ the following inequality holds:
$\|\dot{c}^{(k)}\|_{m-2,\nu}<\varepsilon$. This inequality is
proved by the method of mathematical induction.

Therefore the sequence $\{\dot{c}^{(k)}\}$ is bounded in the space
$C^{m-2,\nu}(\bar{D})$.

We will show that the sequence $\{\dot{c}^{(k)}\}$ is convergent
in the space  $C^{m-2,\nu}(\bar{D})$. Consider the equations:
$$\dot{c}^{(k)}=L_{a}\dot{c}^{(k)} +\gamma_{t}, \eqno(5.15) $$
$$\dot{c}^{(k+1)}=L_{a}\dot{c}^{(k+1)}+\gamma_{t}. \eqno(5.16) $$
Subtracting from the second equation the first one we obtain the
equation:
$$\dot{c}^{(k+1)}-\dot{c}^{(k)}=L_{a}(\dot{c}^{(k+1)})-L_{a}(\dot{c}^{(k)}).
\eqno(5.17) $$

Using the explicit form of $L_{a}$ we have the estimate:
$$\|\dot{c}^{(k+1)}-\dot{c}^{(k)}\|_{m-2,\nu}
\leq K_{3}(t)\|\dot{c}^{(k)}-\dot{c}^{(k-1)}\|_{m-2,\nu},
\eqno(5.18) $$ where we can choose $t_{0}$ such that the following
condition holds $K_{3}(t)<1$ for all $t\in[0,t_{0})$. Then the
sequence $\{\dot{c}^{(k)}\}$ is Cauchy sequence in the space
$C^{m-2,\nu}(\bar{D})$ and therefore is convergent since the space
$C^{m-2,\nu}(\bar{D})$ is complete.

We will show that obtained solution is continuous by $t$. We have:
$$\dot{c}(t_{1})-\dot{c}(t_{2})=L_{a}(\dot{c}(t_{1}))-L_{a}(\dot{c}(t_{2}))+
\gamma_{t_{1}}-\gamma_{t_{2}}. \eqno(5.19) $$

Then there is the estimate:
$$\|\dot{c}(t_{1})-\dot{c}(t_{2})\|_{m-2,\nu}\leq \delta_{1}(t_{1},t_{2})+
\delta_{2}(t_{1},t_{2})\|\dot{c}(t_{1})-\dot{c}(t_{2})\|_{m-2,\nu},
\eqno(5.20) $$ where function $\delta_{1}$ converges to zero if
$|t_{1}-t_{2}|$ converges to zero. Function
$\delta_{2}(t_{1},t_{2})$ is such that
 for any $N>0$ we can choose such
$t_{0}>0$ that for any $t_{1}$ and $t_{2}\in [0,t_{0})$ the
following inequality holds $|\delta_{2}(t_{1},t_{2})|<N$.
Therefore we obtain the continuity of solution.

We will show that the equation $\dot{c}=L_{a}\dot{c} +\gamma_{t}$
has unique solution of class $C^{m-2,\nu}(\bar{D})$ for all
sufficiently small $t\geq 0.$ Let there exist two different
solutions $\dot{c}_{(1)}, \dot{c}_{(2)}$ of class
$C^{m-2,\nu}(\bar{D}).$

Consider the equations:
$$\dot{c}_{(1)}=L_{a}\dot{c}_{(1)} +\gamma_{t}, \eqno(5.21) $$
$$\dot{c}_{(2)}=L_{a}\dot{c}_{(2)}+\gamma_{t}. \eqno(5.22) $$
Subtracting from the second equation the first one we obtain the
equation:
$$\dot{c}_{(2)}-\dot{c}_{(1)}=L_{a}(\dot{c}_{(2)})-L_{a}(\dot{c}_{(1)}).
\eqno(5.23) $$

Using the explicit form of $L_{a}$ we have the estimate:
$$\|\dot{c}_{(2)}-\dot{c}_{(1)}\|_{m-2,\nu}\leq
K_{17}(t)\|\dot{c}_{(2)}-\dot{c}_{(1)}\|_{m-2,\nu}. \eqno(5.24) $$
We can choose $t_{0}$ such that the following condition holds
$K_{17}(t)<1$ for all $t\in[0,t_{0}).$ Therefore we have
contradiction. Therefore $\dot{c}_{(1)}\equiv \dot{c}_{(2)}$ for
all sufficiently small $t\geq 0.$

Lemma 5.1. is proved.

Since curve $l^{*}$ is arbitrary admissible in $D$ therefore
 the equation (5.4) is solvable uniquely
 for any continuous functions $\dot{a}^{i}$ and $\partial_{p}\dot{a}^{i}.$

{\bf Corollary.}
{ \it Let the conditions of lemma 5.1. hold.\\
Then the function $\dot{c}$ takes the form:
$$ \dot{c}(x^{1},x^{2},t)=
\int\limits_{(x^{1}_{(0)},x^{2}_{(0)})}^{(x^{1},x^{2})}
\Biggl(-V\dot{a}^{1}\Biggr) d\tilde{x}^{1}+
\Biggl(-V\dot{a}^{2}\Biggr)d\tilde{x}^{2}+P(\dot{a}^{1},\dot{a}^{2}),
\eqno(5.25) $$ and for $P$ the following inequality holds:
$$\|P(\dot{a}^{1}_{(1)},\dot{a}^{2}_{(1)})-
P(\dot{a}^{1}_{(2)},\dot{a}^{2}_{(2)})\|_{m-2,\nu}\leq
K_{8}(t)(\|\dot{a}^{1}_{(1)}-\dot{a}^{1}_{(2)}\|_{m-2,\nu}+
\|\dot{a}^{2}_{(1)}-\dot{a}^{2}_{(2)}\|_{m-2,\nu}),$$ where for
any $\varepsilon>0$ there exists $t_{0}>0$ such that
 for all $t\in[0,t_{0})$ the following inequality holds:
$K_{8}(t)<\varepsilon.$ }

The proof follows from construction of function $\dot{c}.$

\section*{\S 6. Deduction the formulas of deformations preserving
the product of principal curvatures.}

\section*{\S 6.1. Deduction the formula of $\Delta(g).$}

Consider the following
$$\Delta(g)=g_{t}-g, \eqno(6.1.1)$$
where $g_{t}$ is determinant of the first fundamental form matrix
of hypersurface $F_{t}$.

We will calculate $\Delta(g_{ij}).$ Deformation $\{F_{t}\}$ of
surface $F$ is defined by the formula (1.1). We will use (2.1),
(2.2), (2.3), where
$$a^{j}(0)\equiv 0, c(0)\equiv 0. \eqno(6.1.2)$$
Notice that deformation of surface $F$ determines by the functions
$a^{j}$ and $c.$

Let $\tilde{a}_{\alpha\beta}(t)$ be metric tensor of Riemannian
space at the point $(y^{\sigma}+z^{\sigma}(t)),$
$\tilde{a}_{\alpha\beta}(t)\equiv
\tilde{a}_{\alpha\beta}(y^{\sigma}+z^{\sigma}(t)),$
$\tilde{a}_{\alpha\beta}(0)\equiv
\tilde{a}_{\alpha\beta}(y^{\sigma}).$
$\tilde{a}_{\alpha\beta}\equiv \tilde{a}_{\alpha\beta}(0).$ The
designations $\Gamma_{\beta \sigma}^{ \gamma}(0)$ and
$\Gamma_{\alpha \sigma, \beta}(0)$ mean that the Christoffel
symbols are calculated at the point $(y^{\sigma}).$

$$\Delta(g_{ij})=\tilde{a}_{\alpha\beta}(t)(y^{\alpha},_{i}+z^{\alpha},_{i})
(y^{\beta},_{j}+z^{\beta},_{j})-
\tilde{a}_{\alpha\beta}(0)y^{\alpha},_{i}y^{\beta},_{j}.
\eqno(6.1.3) $$

Then we obtain:\\
$$\Delta(g_{ij})=
(\tilde{a}_{\alpha\beta}(0)+ \frac{\partial
\tilde{a}_{\alpha\beta}}{\partial y^{\sigma}}(0) z^{\sigma})
(y^{\alpha},_{i}+z^{\alpha},_{i}) (y^{\beta},_{j}+z^{\beta},_{j})-
\tilde{a}_{\alpha\beta}(0)y^{\alpha},_{i}y^{\beta},_{j}+$$
$$
 (\tilde{a}_{\alpha\beta}(t)-\tilde{a}_{\alpha\beta}(0)-
\frac{\partial \tilde{a}_{\alpha\beta}}{\partial y^{\sigma}}(0)
z^{\sigma}) (y^{\alpha},_{i}+z^{\alpha},_{i})
(y^{\beta},_{j}+z^{\beta},_{j}). \eqno(6.1.4) $$

Therefore we have:\\
$$\Delta(g_{ij})=
\tilde{a}_{\alpha\beta}(0) (y^{\alpha},_{i}+z^{\alpha},_{i})
(y^{\beta},_{j}+z^{\beta},_{j})-
\tilde{a}_{\alpha\beta}(0)y^{\alpha},_{i}y^{\beta},_{j}+$$
$$
 (\tilde{a}_{\alpha\beta}(t)-\tilde{a}_{\alpha\beta}(0))
(y^{\alpha},_{i}+z^{\alpha},_{i}) (y^{\beta},_{j}+z^{\beta},_{j})+
\frac{\partial \tilde{a}_{\alpha\beta}}{\partial y^{\sigma}}(0)
z^{\sigma} y^{\alpha},_{i}y^{\beta},_{j}-$$
$$\frac{\partial \tilde{a}_{\alpha\beta}}{\partial y^{\sigma}}(0) z^{\sigma}
y^{\alpha},_{i}y^{\beta},_{j}. \eqno(6.1.5) $$

Hence we have:\\
$$\Delta(g_{ij})=
\tilde{a}_{\alpha\beta}(0) (y^{\alpha},_{i}z^{\beta},_{j}+
y^{\beta},_{j}z^{\alpha},_{i})+ \frac{\partial
\tilde{a}_{\alpha\beta}}{\partial y^{\sigma}}(0) z^{\sigma}
y^{\alpha},_{i}y^{\beta},_{j}+
\tilde{a}_{\alpha\beta}(0)z^{\alpha},_{i}z^{\beta},_{j}+$$
$$ (\tilde{a}_{\alpha\beta}(t)-\tilde{a}_{\alpha\beta}(0))
(y^{\alpha},_{i}+z^{\alpha},_{i}) (y^{\beta},_{j}+z^{\beta},_{j})-
\frac{\partial \tilde{a}_{\alpha\beta}}{\partial y^{\sigma}}(0)
z^{\sigma} y^{\alpha},_{i}y^{\beta},_{j}. \eqno(6.1.6)$$

Consider the formula:\\
$$\frac{\partial \tilde{a}_{\alpha\beta}}{\partial y^{\sigma}}(0)=
\Gamma_{\alpha \sigma, \beta}(0)+\Gamma_{\beta \sigma, \alpha}(0)=
\tilde{a}_{\gamma\beta} \Gamma_{\alpha \sigma}^{ \gamma}(0)+
\tilde{a}_{\gamma\alpha} \Gamma_{\beta \sigma}^{ \gamma}(0).
\eqno(6.1.7) $$ where $\Gamma_{\alpha \sigma, \beta}(0) ,
 \Gamma_{\beta \sigma}^{ \gamma}(0)$ are calculated at the point
  $(y^{\sigma}).$

Then we obtain:\\
$$\frac{\partial \tilde{a}_{\alpha\beta}}{\partial y^{\sigma}}(0) z^{\sigma}
y^{\alpha},_{i}y^{\beta},_{j}= \tilde{a}_{\gamma\beta}
\Gamma_{\alpha \sigma}^{ \gamma}(0)
z^{\sigma}y^{\alpha},_{i}y^{\beta},_{j}+ \tilde{a}_{\gamma\alpha}
\Gamma_{\beta \sigma}^{ \gamma}(0)
z^{\sigma}y^{\alpha},_{i}y^{\beta},_{j}. \eqno(6.1.8) $$

We change the positions of indices $\alpha$ and $\gamma$ in the
first term in the right part of the equation (6.1.8) and we also
change the positions of indices $\beta$ and $\gamma$ in the second
term.

Therefore we have:\\
$$\frac{\partial \tilde{a}_{\alpha\beta}}{\partial y^{\sigma}}(0) z^{\sigma}
y^{\alpha},_{i}y^{\beta},_{j}= \tilde{a}_{\alpha\beta}
\Gamma_{\gamma \sigma}^{\alpha}(0)
z^{\sigma}y^{\gamma},_{i}y^{\beta},_{j}+ \tilde{a}_{\alpha\beta}
\Gamma_{\gamma \sigma}^{ \beta}(0)
z^{\sigma}y^{\alpha},_{i}y^{\gamma},_{j}. \eqno(6.1.9)$$

Considering the following formula (2.2) we have:\\
$$\Delta(g_{ij})=
\tilde{a}_{\alpha\beta}(0)y^{\alpha},_{i}\nabla^{*}_{j}z^{\beta}+
\tilde{a}_{\alpha\beta}(0)y^{\beta},_{j}\nabla^{*}_{i}z^{\alpha}+
 (\tilde{a}_{\alpha\beta}(t)-\tilde{a}_{\alpha\beta}(0)-
\frac{\partial \tilde{a}_{\alpha\beta}}{\partial y^{\sigma}}(0)
z^{\sigma}) y^{\alpha},_{i}y^{\beta},_{j}+$$
$$
 (\tilde{a}_{\alpha\beta}(t)-\tilde{a}_{\alpha\beta}(0))
(y^{\alpha},_{i}z^{\beta},_{j}+y^{\beta},_{j}z^{\alpha},_{i})+
\tilde{a}_{\alpha\beta}(t)z^{\alpha},_{i}z^{\beta},_{j}.
\eqno(6.1.10) $$

Then we obtain:\\
$$g^{ij}\Delta(g_{ij})=
2g^{ij}\tilde{a}_{\alpha\beta}(0)y^{\alpha},_{i}\nabla^{*}_{j}z^{\beta}+
 (\tilde{a}_{\alpha\beta}(t)-\tilde{a}_{\alpha\beta}(0)-
\frac{\partial \tilde{a}_{\alpha\beta}}{\partial y^{\sigma}}(0)
z^{\sigma}) g^{ij}y^{\alpha},_{i}y^{\beta},_{j}+$$
$$
 2(\tilde{a}_{\alpha\beta}(t)-\tilde{a}_{\alpha\beta}(0))
g^{ij}y^{\alpha},_{i}z^{\beta},_{j}+
\tilde{a}_{\alpha\beta}(t)g^{ij}z^{\alpha},_{i}z^{\beta},_{j}.
\eqno(6.1.11) $$

Denote: \\
$$W_{1}= (\tilde{a}_{\alpha\beta}(t)-\tilde{a}_{\alpha\beta}(0)-
\frac{\partial \tilde{a}_{\alpha\beta}}{\partial y^{\sigma}}(0)
z^{\sigma}) g^{ij}y^{\alpha},_{i}y^{\beta},_{j}+$$
$$
 2(\tilde{a}_{\alpha\beta}(t)-\tilde{a}_{\alpha\beta}(0))
g^{ij}y^{\alpha},_{i}z^{\beta},_{j}+
\tilde{a}_{\alpha\beta}(t)g^{ij}z^{\alpha},_{i}z^{\beta},_{j}.
\eqno(6.1.12) $$

Using the properties of determinant we have:
$$\Delta(g)=gg^{ij}\Delta(g_{ij})+W_{2}. \eqno(6.1.13) $$
where $$W_{2}=\Delta(g_{11}) \Delta(g_{22})-(\Delta(g_{12}))^2.
\eqno(6.1.14)$$

Then the equation (6.1.2) takes the form:
$$\Delta(g)=2gg^{ij}\tilde{a}_{\alpha\beta}y^{\alpha},_{i}\nabla^{*}_{j}z^{\beta}+
gW_{1}+W_{2}. \eqno(6.1.15) $$

Using the equation (6.1.3) we write the equation (6.1.15) as:
$$\frac{\Delta(g)}{2g}=a^{l},_{l}-cb_{lm}g^{ml}+\frac{W_{1}}{2}+
\frac{W_{2}}{2g}. \eqno(6.1.16)$$ Using the formula
$\partial_{i}(\ln\sqrt{g})=\Gamma^{j}_{ij},$ where
$\Gamma^{k}_{ij}$ are the Christoffel symbols for hypersurface $F$
in the metric $g_{ij}$ and formula of mean curvature
$2H=g^{im}b_{im}$ we write the equation (6.1.16) as
$$\frac{\Delta(g)}{2\sqrt{g}}=
\partial_{l}(\sqrt{g}a^{l})-2Hc\sqrt{g}+
\frac{\sqrt{g}W_{1}}{2}+\frac{W_{2}}{2\sqrt{g}}. \eqno(6.1.17) $$
The equation (6.1.17) is required equation for functions $a^{i}$
and $c$, determining continuous $A-$deformation of hypersurface
$F$.

Equation (6.1.17) takes the form:
$$\frac{\Delta(g)}{2g}=
 \partial_{1}a^{1}+\partial_{2}a^{2}+a^{1}\partial_{1}(\ln\sqrt{g})+
a^{2}\partial_{2}(\ln\sqrt{g})-\Psi_{2} , \eqno(6.1.18) $$ where
$$\Psi_{2}=2Hc-\frac{W_{1}}{2}-\frac{W_{2}}{2g}. \eqno(6.1.19)$$

Then we obtain:
$$\Delta(g)=2g(
\partial_{1}a^{1}+\partial_{2}a^{2}+q_{k}a^{k}-\Psi_{2}), \eqno(6.1.20)$$
where $$q_{1}=\partial_{1}(\ln \sqrt{g}), q_{2}=\partial_{2}(\ln
\sqrt{g}). \eqno(6.1.21)$$ Note that $q_{k}$ do not depend on $t.$

Equation (6.1.20) determines $\Delta(g)$ for deformations of
surface $F$ in $R^3.$

\section*{\S 6.2. Deduction the formulas of deformations preserving
the product of principal curvatures.}

Deformation $\{F_{t}\}$ of surface $F$ is determined by (2.1). We
will deduct the formulas of changing
 of the second fundamental form determinant.

The condition of preservation the product of principal curvatures
takes the following form:
$$\Delta(g)=\frac{g}{b}\Delta(b). \eqno(6.2.1) $$
$$\Delta(b)=bb^{ij}\Delta(b_{ij})+W_{2}^{(b)}. \eqno(6.2.2) $$
We have the formula:\\
$$\Delta(K)=\frac{1}{b(t)}(\Delta(g)-\frac{g}{b}\Delta(b)),
b(t)=b+\Delta(b). \eqno(6.2.3) $$

We introduce conjugate isothermal coordinate system where \\
$$b_{ii}=V,i=1,2, b_{12}=b_{21}=0,
b^{ii}=\frac{1}{V},i=1,2, b^{12}=b^{21}=0. \eqno(6.2.4) $$

Then we have:
$$\Delta(b)=V(\Delta(b_{11})+\Delta(b_{22}))+W_{2}^{(b)}, \eqno(6.2.5) $$
where
$$W_{2}^{(b)}=\Delta(b_{11}) \Delta(b_{22})-(\Delta(b_{12}))^2. \eqno(6.2.6)$$

Therefore the condition of preservation the product of principal
curvatures takes the following form:
$$\Delta(g)=\frac{g}{V}(\Delta(b_{11})+\Delta(b_{22}))+
\frac{g}{V^{2}}W_{2}^{(b)}. \eqno(6.2.7) $$

We have the following formula:
$$b_{ij}(0)=-\tilde{a}_{\alpha\beta}(0)y^{\alpha}_{,i} \nabla^{*}_{j}n^{\beta}(0). \eqno(6.2.8) $$

$$b_{ij}(t)=
-\tilde{a}_{\alpha\beta}(t)(y^{\alpha}_{,i}+z^{\alpha}_{,i})
\nabla^{*}_{j}\tilde{n}^{\beta}(t), \eqno(6.2.9) $$ where
$\tilde{n}^{\beta}(t)$ is unit normal vector at the point
$(y^{\alpha}+z^{\alpha}).$

Then we obtain:
$$b_{ij}(t)=
-\tilde{a}_{\alpha\beta}(t)(y^{\alpha}_{,i}+z^{\alpha}_{,i})
(\tilde{n}_{,j}^{\beta}(t)+
\Gamma^{\beta}_{\mu\sigma}(t)(y^{\mu}_{,j}+z^{\mu}_{,j})
\tilde{n}^{\sigma}(t)). \eqno(6.2.10) $$

Let $n^{\beta}(t)$ be result of parallel transfer of unite normal
vector
 $n^{\beta}(0)$ to the point $(y^{\alpha}+z^{\alpha})$ along
the path of translation by deformation. Therefore we have the
following formula for all sufficiently small $t$:
$$n^{\alpha}(t)=n^{\alpha}(0)+
n^{\sigma}(0)\sum_{k=1}^{\infty}A^{\alpha}_{(k)\sigma}(0,t).
\eqno(6.2.11) $$

Use the following formula:\\
$$\tilde{n}^{\beta}(t)=
\frac{n^{\beta}(t)}
{\sqrt{\tilde{a}_{\alpha_{0}\beta_{0}}(t)n^{\alpha_{0}}(t)n^{\beta_{0}}(t)}}.
\eqno(6.2.12) $$
Denote:\\
$$\|n(t)\|=
\sqrt{\tilde{a}_{\alpha_{0}\beta_{0}}(t)n^{\alpha_{0}}(t)n^{\beta_{0}}(t)}.
\eqno(6.2.13) $$

Then we have:
$$b_{ij}(t)=
-\tilde{a}_{\alpha\beta}(t)z^{\alpha}_{,i}\nabla^{*}_{j}\tilde{n}^{\beta}(t)
-\tilde{a}_{\alpha\beta}(t)y^{\alpha}_{,i}\nabla^{*}_{j}\tilde{n}^{\beta}(t).
\eqno(6.2.14) $$ Using the formulas (2.27) and (2.31) we obtain:
$$-\tilde{a}_{\alpha\beta}(t)z^{\alpha}_{,i}\nabla^{*}_{j}\tilde{n}^{\beta}(t)=
-\tilde{a}_{\alpha\beta}(t)\partial_{i}(a^{k})y^{\alpha},_{k}
\nabla^{*}_{j}\tilde{n}^{\beta}(t)+M^{1}_{ij}, \eqno(6.2.15) $$
where $$M^{1}_{ij}=-\tilde{a}_{\alpha\beta}(t)
(\Gamma^{k}_{pi}a^{p}y^{\alpha},_{k}+c,_{i}n^{\alpha}+
a^{k}y^{\alpha}_{,k,i}+
cn^{\alpha}_{,i})\nabla^{*}_{j}\tilde{n}^{\beta}(t).
\eqno(6.2.16)$$

Then we have:
$$-\tilde{a}_{\alpha\beta}(t)z^{\alpha}_{,i}\nabla^{*}_{j}\tilde{n}^{\beta}(t)=
-\tilde{a}_{\alpha\beta}(0)\partial_{i}(a^{k})y^{\alpha},_{k}
\nabla^{*}_{j} n^{\beta}(0)-$$
$$\tilde{a}_{\alpha\beta}(t)\partial_{i}(a^{k})y^{\alpha},_{k}
\nabla^{*}_{j}\tilde{n}^{\beta}(t)+
\tilde{a}_{\alpha\beta}(0)\partial_{i}(a^{k})y^{\alpha},_{k}
\nabla^{*}_{j} n^{\beta}(0)+M^{1}_{ij}, \eqno(6.2.17) $$

Define:\\
$$ M^{2}_{ij}=-\tilde{a}_{\alpha\beta}(t)\partial_{i}(a^{k})y^{\alpha},_{k}
\nabla^{*}_{j}\tilde{n}^{\beta}(t)+
\tilde{a}_{\alpha\beta}(0)\partial_{i}(a^{k})y^{\alpha},_{k}
\nabla^{*}_{j} n^{\beta}(0)+M^{1}_{ij}. \eqno(6.2.18) $$

Consequently we get:
$$-\tilde{a}_{\alpha\beta}(t)z^{\alpha}_{,i}\nabla^{*}_{j}\tilde{n}^{\beta}(t)=
-\tilde{a}_{\alpha\beta}(0)\partial_{i}(a^{k})y^{\alpha},_{k}
\nabla^{*}_{j} n^{\beta}(0)+M^{2}_{ij}. \eqno(6.2.19) $$

We use the following equation:
$$ \nabla^{*}_{j} n^{\beta}(0)=-b_{jk}g^{kl}y^{\beta}_{,l}. \eqno(6.2.20)$$

Then we obtain:
$$-\tilde{a}_{\alpha\beta}(t)z^{\alpha}_{,i}\nabla^{*}_{j}\tilde{n}^{\beta}(t)=
\partial_{i}(a^{k})b_{jk}+M^{2}_{ij}. \eqno(6.2.21) $$

Using the fact $b_{12}=0$ we have:\\
$$-\tilde{a}_{\alpha\beta}(t)z^{\alpha}_{,1}\nabla^{*}_{1}\tilde{n}^{\beta}(t)=
V \partial_{1}(a^{1})+M^{2}_{11}, \eqno(6.2.22) $$
$$-\tilde{a}_{\alpha\beta}(t)z^{\alpha}_{,2}\nabla^{*}_{2}\tilde{n}^{\beta}(t)=
V \partial_{2}(a^{2})+M^{2}_{22}. \eqno(6.2.23) $$

We have the expression:\\
$$\nabla^{*}_{j}\tilde{n}^{\beta}(t)=
\tilde{n}_{,j}^{\beta}(t)+
\Gamma^{\beta}_{\mu\sigma}(t)(y^{\mu}_{,j}+z^{\mu}_{,j})
\tilde{n}^{\sigma}(t)=
\Biggl(\frac{n^{\beta}(t)}{\|n(t)\|}\Biggr)_{,j}+
\Gamma^{\beta}_{\mu\sigma}(t)(y^{\mu}_{,j}+z^{\mu}_{,j})
\Biggl(\frac{n^{\sigma}(t)}{\|n(t)\|}\Biggr). \eqno(6.2.24) $$

Then we obtain:\\
$$\nabla^{*}_{j}\tilde{n}^{\beta}(t)=
\frac{\nabla^{*}_{j}n^{\beta}(t)}{\|n(t)\|}+
n^{\beta}(t)\Biggl(\frac{1}{\|n(t)\|}\Biggr)_{,j}. \eqno(6.2.25)$$

Consider the equation:\\
$$-\tilde{a}_{\alpha\beta}(t)y^{\alpha}_{,i}\nabla^{*}_{j}\tilde{n}^{\beta}(t)=
-\tilde{a}_{\alpha\beta}(t)y^{\alpha}_{,i}
\nabla^{*}_{j}n^{\beta}(t)\Biggl(\frac{1}{\|n(t)\|}\Biggr)-
\tilde{a}_{\alpha\beta}(t)y^{\alpha}_{,i}
n^{\beta}(t)\Biggl(\frac{1}{\|n(t)\|}\Biggr)_{,j}. \eqno(6.2.26)$$

Consider the formula:\\
$$\nabla^{*}_{j}n^{\beta}(t)=
n_{,j}^{\beta}(t)+
\Gamma^{\beta}_{\mu\sigma}(t)(y^{\mu}_{,j}+z^{\mu}_{,j})n^{\sigma}(t).
\eqno(6.2.27) $$

Then we have:\\
$$\nabla^{*}_{j}n^{\beta}(t)=
n_{,j}^{\beta}(t)+
\Gamma^{\beta}_{\mu\sigma}(t)y^{\mu}_{,j}n^{\sigma}(t)+
\Gamma^{\beta}_{\mu\sigma}(t)z^{\mu}_{,j}n^{\sigma}(t)=$$
$$=n_{,j}^{\beta}(t)+
\Gamma^{\beta}_{\mu\sigma}(0)y^{\mu}_{,j}n^{\sigma}(t)+
(\Gamma^{\beta}_{\mu\sigma}(t)-\Gamma^{\beta}_{\mu\sigma}(0))
y^{\mu}_{,j}n^{\sigma}(t)+$$
$$\Gamma^{\beta}_{\mu\sigma}(0)z^{\mu}_{,j}n^{\sigma}(t)+
(\Gamma^{\beta}_{\mu\sigma}(t)-\Gamma^{\beta}_{\mu\sigma}(0))
z^{\mu}_{,j}n^{\sigma}(t). \eqno(6.2.28) $$

We will use the following formula:\\
$$n^{\beta}(t)=n^{\beta}(0)+
n^{\sigma}(0)\sum_{k=1}^{\infty}A^{\beta}_{(k)\sigma}(0,t).
\eqno(6.2.29) $$
Denote:\\
$$A_{1}^{\beta}(t)=
n^{\sigma}(0)\sum_{k=1}^{\infty}A^{\beta}_{(k)\sigma}(0,t),
\eqno(6.2.30) $$
$$A_{2}^{\beta}(t)=
n^{\sigma}(0)\sum_{k=2}^{\infty}A^{\beta}_{(k)\sigma}(0,t).
\eqno(6.2.31) $$

We use the equation:\\
$$n^{\beta}(t)=n^{\beta}(0)+
n^{\sigma}(0)A^{\beta}_{(1)\sigma}(0,t)+A_{2}^{\beta}(t).
\eqno(6.2.32) $$

We have:\\
$$n^{\beta}_{,j}(t)=n^{\beta}_{,j}(0)+
n^{\sigma}_{,j}(0)A^{\beta}_{(1)\sigma}(0,t)+
n^{\sigma}(0)A^{\beta}_{(1)\sigma,j}(0,t)+A_{2,j}^{\beta}(t).
\eqno(6.2.33) $$

Hence:\\
$$\nabla^{*}_{j}n^{\beta}(t)=
n^{\beta}_{,j}(0)+ n^{\sigma}_{,j}(0)A^{\beta}_{(1)\sigma}(0,t)+
n^{\sigma}(0)A^{\beta}_{(1)\sigma,j}(0,t)+A_{2,j}^{\beta}(t)+$$
$$\Gamma^{\beta}_{\mu\tau}(0)y^{\mu}_{,j}n^{\tau}(0)+
\Gamma^{\beta}_{\mu\tau}(0)y^{\mu}_{,j}
n^{\sigma}(0)A^{\tau}_{(1)\sigma}(0,t)+
\Gamma^{\beta}_{\mu\tau}(0)y^{\mu}_{,j}A_{2}^{\tau}(t)+$$
$$
(\Gamma^{\beta}_{\mu\tau}(t)-\Gamma^{\beta}_{\mu\tau}(0))
y^{\mu}_{,j}n^{\tau}(0)+
(\Gamma^{\beta}_{\mu\tau}(t)-\Gamma^{\beta}_{\mu\tau}(0))
y^{\mu}_{,j}A_{1}^{\tau}(t)+$$
$$\Gamma^{\beta}_{\mu\tau}(0)z^{\mu}_{,j}n^{\tau}(0)+
\Gamma^{\beta}_{\mu\tau}(0)z^{\mu}_{,j}A_{1}^{\tau}(t)+
(\Gamma^{\beta}_{\mu\tau}(t)-\Gamma^{\beta}_{\mu\tau}(0))
z^{\mu}_{,j}n^{\tau}(t). \eqno(6.2.34) $$

Denote:\\
$$T^{\beta}_{j}=
n^{\sigma}_{,j}(0)A^{\beta}_{(1)\sigma}(0,t)+ A_{2,j}^{\beta}(t)+
\Gamma^{\beta}_{\mu\tau}(0)y^{\mu}_{,j}
n^{\sigma}(0)A^{\tau}_{(1)\sigma}(0,t)+
\Gamma^{\beta}_{\mu\tau}(0)y^{\mu}_{,j}A_{2}^{\tau}(t)+$$
$$
(\Gamma^{\beta}_{\mu\tau}(t)-\Gamma^{\beta}_{\mu\tau}(0))
y^{\mu}_{,j}n^{\tau}(0)+
(\Gamma^{\beta}_{\mu\tau}(t)-\Gamma^{\beta}_{\mu\tau}(0))
y^{\mu}_{,j}A_{1}^{\tau}(t)+$$
$$\Gamma^{\beta}_{\mu\tau}(0)z^{\mu}_{,j}A_{1}^{\tau}(t)+
(\Gamma^{\beta}_{\mu\tau}(t)-\Gamma^{\beta}_{\mu\tau}(0))
z^{\mu}_{,j}n^{\tau}(t). \eqno(6.2.35) $$

Then we get:\\
$$\nabla^{*}_{j}n^{\beta}(t)=
n^{\beta}_{,j}(0)+
\Gamma^{\beta}_{\mu\tau}(0)y^{\mu}_{,j}n^{\tau}(0)+
n^{\sigma}(0)A^{\beta}_{(1)\sigma,j}(0,t)+
\Gamma^{\beta}_{\mu\tau}(0)z^{\mu}_{,j}n^{\tau}(0)+T^{\beta}_{j}=$$
$$\nabla^{*}_{j}n^{\beta}(0)+
n^{\sigma}(0)A^{\beta}_{(1)\sigma,j}(0,t)+
\Gamma^{\beta}_{\mu\tau}(0)z^{\mu}_{,j}n^{\tau}(0)+T^{\beta}_{j}.
\eqno(6.2.36) $$

Consider the expression:\\
$$-\tilde{a}_{\alpha\beta}(t)y^{\alpha}_{,i}\nabla^{*}_{j}n^{\beta}(t)=
-\tilde{a}_{\alpha\beta}(0)y^{\alpha}_{,i}\nabla^{*}_{j}n^{\beta}(t)-
(\tilde{a}_{\alpha\beta}(t)-\tilde{a}_{\alpha\beta}(0))
y^{\alpha}_{,i}\nabla^{*}_{j}n^{\beta}(t)=$$
$$\tilde{a}_{\alpha\beta}(0)y^{\alpha}_{,i}b_{jk}g^{kl}y^{\beta}_{,l}-
\tilde{a}_{\alpha\beta}(0)y^{\alpha}_{,i}
(n^{\sigma}(0)A^{\beta}_{(1)\sigma,j}(0,t)+
\Gamma^{\beta}_{\mu\tau}(0)z^{\mu}_{,j}n^{\tau}(0))-$$
$$\tilde{a}_{\alpha\beta}(0)y^{\alpha}_{,i}T^{\beta}_{j}-
(\tilde{a}_{\alpha\beta}(t)-\tilde{a}_{\alpha\beta}(0))
y^{\alpha}_{,i}\nabla^{*}_{j}n^{\beta}(t). \eqno(6.2.37) $$

Then we obtain:\\
$$-\tilde{a}_{\alpha\beta}(t)y^{\alpha}_{,i}\nabla^{*}_{j}n^{\beta}(t)=
b_{ji}- \tilde{a}_{\alpha\beta}(0)y^{\alpha}_{,i}
(n^{\sigma}(0)A^{\beta}_{(1)\sigma,j}(0,t)+
\Gamma^{\beta}_{\mu\tau}(0)z^{\mu}_{,j}n^{\tau}(0))-$$
$$\tilde{a}_{\alpha\beta}(0)y^{\alpha}_{,i}T^{\beta}_{j}-
(\tilde{a}_{\alpha\beta}(t)-\tilde{a}_{\alpha\beta}(0))
y^{\alpha}_{,i}\nabla^{*}_{j}n^{\beta}(t). \eqno(6.2.38) $$

Therefore:\\
$$-\tilde{a}_{\alpha\beta}(t)y^{\alpha}_{,i}\nabla^{*}_{j}n^{\beta}(t)
\Biggl(\frac{1}{\|n(t)\|}\Biggr)=$$
$$b_{ji}\Biggl(\frac{1}{\|n(t)\|}\Biggr)-
\tilde{a}_{\alpha\beta}(0)y^{\alpha}_{,i}
(n^{\sigma}(0)A^{\beta}_{(1)\sigma,j}(0,t)+
\Gamma^{\beta}_{\mu\tau}(0)z^{\mu}_{,j}n^{\tau}(0))
\Biggl(\frac{1}{\|n(t)\|}\Biggr)-$$
$$\tilde{a}_{\alpha\beta}(0)y^{\alpha}_{,i}T^{\beta}_{j}
\Biggl(\frac{1}{\|n(t)\|}\Biggr)-
(\tilde{a}_{\alpha\beta}(t)-\tilde{a}_{\alpha\beta}(0))
y^{\alpha}_{,i}\nabla^{*}_{j}n^{\beta}(t)\Biggl(\frac{1}{\|n(t)\|}\Biggr).
\eqno(6.2.39) $$

We change the form of last expression:\\
$$-\tilde{a}_{\alpha\beta}(t)y^{\alpha}_{,i}\nabla^{*}_{j}n^{\beta}(t)
\Biggl(\frac{1}{\|n(t)\|}\Biggr)=$$
$$b_{ji}+
\frac{b_{ji}(1-\|n(t)\|)}{\|n(t)\|}-
\tilde{a}_{\alpha\beta}(0)y^{\alpha}_{,i}
(n^{\sigma}(0)A^{\beta}_{(1)\sigma,j}(0,t)+
\Gamma^{\beta}_{\mu\tau}(0)z^{\mu}_{,j}n^{\tau}(0))
\Biggl(\frac{1}{\|n(t)\|}\Biggr)-$$
$$\tilde{a}_{\alpha\beta}(0)y^{\alpha}_{,i}T^{\beta}_{j}
\Biggl(\frac{1}{\|n(t)\|}\Biggr)-
(\tilde{a}_{\alpha\beta}(t)-\tilde{a}_{\alpha\beta}(0))
y^{\alpha}_{,i}\nabla^{*}_{j}n^{\beta}(t)\Biggl(\frac{1}{\|n(t)\|}\Biggr).
\eqno(6.2.40) $$

Define:\\
$$M^{3}_{ij}=
\frac{b_{ji}(1-\|n(t)\|)}{\|n(t)\|}-
\tilde{a}_{\alpha\beta}(0)y^{\alpha}_{,i}
(n^{\sigma}(0)A^{\beta}_{(1)\sigma,j}(0,t)+
\Gamma^{\beta}_{\mu\tau}(0)z^{\mu}_{,j}n^{\tau}(0))
\Biggl(\frac{1}{\|n(t)\|}\Biggr)-$$
$$\tilde{a}_{\alpha\beta}(0)y^{\alpha}_{,i}T^{\beta}_{j}
\Biggl(\frac{1}{\|n(t)\|}\Biggr)-
(\tilde{a}_{\alpha\beta}(t)-\tilde{a}_{\alpha\beta}(0))
y^{\alpha}_{,i}\nabla^{*}_{j}n^{\beta}(t)\Biggl(\frac{1}{\|n(t)\|}\Biggr).
\eqno(6.2.41) $$

Then we have:\\
$$-\tilde{a}_{\alpha\beta}(t)y^{\alpha}_{,i}\nabla^{*}_{j}n^{\beta}(t)
\Biggl(\frac{1}{\|n(t)\|}\Biggr)=b_{ij}+M^{3}_{ij}. \eqno(6.2.42)
$$

Consequently:\\
$$b_{ij}(t)=
-\tilde{a}_{\alpha\beta}(t)z^{\alpha}_{,i}\nabla^{*}_{j}\tilde{n}^{\beta}(t)-
\tilde{a}_{\alpha\beta}(t)y^{\alpha}_{,i}\nabla^{*}_{j}\tilde{n}^{\beta}(t)=$$
$$\partial_{i}(a^{k})b_{jk}+M^{2}_{ij}+b_{ij}+M^{3}_{ij}-
\tilde{a}_{\alpha\beta}(t)y^{\alpha}_{,i}
n^{\beta}(t)\Biggl(\frac{1}{\|n(t)\|}\Biggr)_{,j}. \eqno(6.2.43)
$$

Denote:\\
$$M^{4}_{ij}=M^{2}_{ij}+M^{3}_{ij}-\tilde{a}_{\alpha\beta}(t)y^{\alpha}_{,i}
n^{\beta}(t)\Biggl(\frac{1}{\|n(t)\|}\Biggr)_{,j}. \eqno(6.2.44)
$$

Hence:\\
$$b_{ij}(t)=\partial_{i}(a^{k})b_{jk}+b_{ij}+M^{4}_{ij}. \eqno(6.2.45) $$

Therefore:\\
$$\Delta(b_{ij})=\partial_{i}(a^{k})b_{jk}+M^{4}_{ij}. \eqno(6.2.46) $$

Then we have:\\
$$\Delta(b_{11})=V\partial_{1}(a^{1})+M^{4}_{11}, \eqno(6.2.47) $$
$$\Delta(b_{22})=V\partial_{2}(a^{2})+M^{4}_{22}. \eqno(6.2.48) $$

Hence the condition of preservation the product of principal
curvatures takes the following form:
$$\Delta(g)=g(\partial_{1}(a^{1})+\partial_{2}(a^{2}))+
\frac{g}{V}(M^{4}_{11}+M^{4}_{22})+\frac{g}{V^{2}}W_{2}^{(b)}.
\eqno(6.2.49) $$

Using the formula (6.1.20)
we obtain the equation of preservation the product of principal curvatures:\\
$$\partial_{1}a^{1}+\partial_{2}a^{2}+2q_{k}a^{k}-2\Psi_{2}=
\frac{1}{V}(M^{4}_{11}+M^{4}_{22})+\frac{1}{V^{2}}W_{2}^{(b)}.
\eqno(6.2.50) $$

Then we have:
$$\partial_{1}a^{1}+\partial_{2}a^{2}+2q_{k}a^{k}=2\Psi_{2}+
\frac{1}{V}(M^{4}_{11}+M^{4}_{22})+\frac{1}{V^{2}}W_{2}^{(b)}.
\eqno(6.2.51) $$

We differentiate the equation (6.2.51) by $t.$ Then we have:
$$\partial_{1}\dot{a}^{1}+\partial_{2}\dot{a}^{2}+2q_{k}\dot{a}^{k}=
2\dot{\Psi}_{2}+\frac{1}{V}(\dot{M}^{4}_{11}+\dot{M}^{4}_{22})+
\frac{1}{V^{2}}\dot{W}_{2}^{(b)}. \eqno(6.2.52) $$

The equation takes the following form:
$$\partial_{1}\dot{a}^{1}+\partial_{2}\dot{a}^{2}+q^{(b)}_{k}\dot{a}^{k}=
\dot{\Psi}_{2}^{(b)},  \eqno(6.2.53)$$ where
$\dot{\Psi}_{2}^{(b)}=q^{(b)}_{0}\dot{c}- P_{0}(\dot{a}^{1},
\dot{a}^{2},\partial_{i}\dot{a}^{j}).$ $P_{0}$ has explicit form.
Notice that $q^{(b)}_{k}\in C^{m-3,\nu},$ $q^{(b)}_{0} \in
C^{m-3,\nu}$ and
 do not depend on $t.$

{\bf Lemma 6.2.1.} {\it Let the following conditions hold:

1) metric tensor in $R^{3}$ satisfies the conditions: $\exists
M_{0}=const>0$ such that
$\|\tilde{a}_{\alpha\beta}\|_{m,\nu}<M_{0},$ $\|\partial
\tilde{a}_{\alpha\beta}\|_{m,\nu}<M_{0},$ $\|\partial^{2}
\tilde{a}_{\alpha\beta}\|_{m,\nu}<M_{0}.$

2) $\exists t_{0}>0$ such that $a^{k}(t), \partial_{i}a^{k}(t),
\dot{a}^{k}(t), \partial_{i}\dot{a}^{k}(t)$ are continuous by $t,
\forall t\in[0,t_{0}],$ $ a^{k}(0)\equiv 0,
\partial_{i}a^{k}(0)\equiv 0.$

3) $\exists t_{0}>0$ such that $a^{i}(t)\in C^{m-2,\nu} ,
\partial_{k} a^{i}(t)\in C^{m-3,\nu},$ $\forall t\in[0,t_{0}].$

Then $\exists t_{*}>0$ such that for all $t\in[0,t_{*})$ $P_{0}\in
C^{m-3,\nu}$ and the following inequality holds:
$$\|P_{0}(\dot{a}^{1}_{(1)},\dot{a}^{2}_{(1)})-
P_{0}(\dot{a}^{1}_{(2)},\dot{a}^{2}_{(2)})\|_{m-2,\nu}\leq
K_{9}(t)(\|\dot{a}^{1}_{(1)}-\dot{a}^{1}_{(2)}\|_{m-1,\nu}+
\|\dot{a}^{2}_{(1)}-\dot{a}^{2}_{(2)}\|_{m-1,\nu}),$$ where for
any $\varepsilon>0$ there exists $t_{0}>0$ such that
 for all $t\in[0,t_{0})$ the following inequality holds:
$K_{9}(t)<\varepsilon.$ }

The proof follows from construction of function $P_{0}$ and lemmas
of \S 7 and \S 8.

The equation (6.2.53) determines deformations of surface $F$
preserving
 the product of principal curvatures with condition of $G-$deformation.

\section*{6.2.1. The formulas of $\Delta(K)$ and $\dot{\Delta}(K).$}
Consider the following formula:
$$\Delta(K)=\frac{1}{b(t)}(\Delta(g)-\frac{g}{b}\Delta(b))=$$
$$\frac{g}{b(t)}(\partial_{1}a^{1}+\partial_{2}a^{2}+2q_{k}a^{k}-(2\Psi_{2}+
\frac{1}{V}(M^{4}_{11}+M^{4}_{22})+\frac{1}{V^{2}}W_{2}^{(b)})),
\eqno(6.2.54) $$
$$b(t)=b+\Delta(b). \eqno(6.2.55) $$

We have:
$$K(t)=K+\Delta(K). \eqno(6.2.56) $$

Therefore: $$\dot{b}(t)=\dot{\Delta}(b),
\dot{K}(t)=\dot{\Delta}(K). \eqno(6.2.57) $$

Then we obtain:
$$\dot{\Delta}(K)=
-\frac{g\dot{\Delta}b(t)}{(b(t))^{2}}
(\partial_{1}a^{1}+\partial_{2}a^{2}+2q_{k}a^{k}-(2\Psi_{2}+
\frac{1}{V}(M^{4}_{11}+M^{4}_{22})+\frac{1}{V^{2}}W_{2}^{(b)}))+$$
$$\frac{g}{b(t)}
(\partial_{1}\dot{a}^{1}+\partial_{2}\dot{a}^{2}+2q_{k}\dot{a}^{k}-
(2\dot{\Psi}_{2}+\frac{1}{V}(\dot{M}^{4}_{11}+\dot{M}^{4}_{22})+
\frac{1}{V^{2}}\dot{W}_{2}^{(b)})). \eqno(6.2.58) $$

We finally obtain the following formula:
$$\dot{\Delta}(K)=
-\frac{g\dot{\Delta}b(t)}{(b(t))^{2}}
(\partial_{1}a^{1}+\partial_{2}a^{2}+2q_{k}a^{k}-(2\Psi_{2}+
\frac{1}{V}(M^{4}_{11}+M^{4}_{22})+\frac{1}{V^{2}}W_{2}^{(b)}))+$$
$$\frac{g}{b(t)}
(\partial_{1}\dot{a}^{1}+\partial_{2}\dot{a}^{2}+q^{(b)}_{k}\dot{a}^{k}
-\dot{\Psi}_{2}^{(b)}). \eqno(6.2.59) $$

\section*{\S 7. Auxiliary estimations of norms.}

Denote:\\
$\|\partial z \|^{(t)}_{m_{1},\nu}= \max_{\alpha,
i}\|z^{\alpha}_{,i} \|^{(t)}_{m_{1},\nu}= \max_{\alpha,
i}\max_{\tau \in [0;t]} \|z^{\alpha}_{,i}(\tau)\|_{m_{1},\nu}.$

{\bf Lemma 7.1.} {it  The following estimations hold:

1) $\| z \|^{(t)}_{m_{1},\nu}\leq M_{5} (\| a
\|^{(t)}_{m_{1},\nu}+\| c \|^{(t)}_{m_{1},\nu}),$

2) $\| a \|^{(t)}_{m_{1},\nu}\leq M_{6}\| z \|^{(t)}_{m_{1},\nu},$

3) $\| c \|^{(t)}_{m_{1},\nu}\leq M_{7}\| z \|^{(t)}_{m_{1},\nu},$

4) $\|\partial z \|^{(t)}_{m_{1},\nu}\leq M_{8}\| z
\|^{(t)}_{m_{1}+1,\nu},$

where constants $M_{5}, M_{6}, M_{7}, M_{8} $ are determined by
surface $F$ and do not depend on  $t.$ }

Proof of lemma follows from properties of norm in the space
$C^{m_{1},\nu}.$

{\bf Lemma 7.2.} {\it The following estimations hold:

1) $ \|W_{1}\|^{(t)}_{m_{1},\nu} \leq M_{2} ((\|z
\|^{(t)}_{m_{1},\nu})^{2}+ \|z \|^{(t)}_{m_{1},\nu}\|\partial z
\|^{(t)}_{m_{1},\nu}+ (\|\partial z \|^{(t)}_{m_{1},\nu})^{2}),$

2) \\
$ \|\Delta(g_{ij})\|^{(t)}_{m_{1},\nu} \leq M_{3} (\|z
\|^{(t)}_{m_{1},\nu}+\|\partial z \|^{(t)}_{m_{1},\nu}+ (\|z
\|^{(t)}_{m_{1},\nu})^{2}+ \|z \|^{(t)}_{m_{1},\nu}\|\partial z
\|^{(t)}_{m_{1},\nu}+ (\|\partial z \|^{(t)}_{m_{1},\nu})^{2}).$

3) $ \|W_{2}\|^{(t)}_{m_{1},\nu} \leq M_{4} (\max_{i,j}
\|\Delta(g_{ij})\|^{(t)}_{m_{1},\nu})^{2},$

where constants $M_{2}, M_{3}, M_{4}$ are determined by surface
$F$ and do not depend on $t.$ }

Proof of lemma follows from properties of norms
 in the space $C^{m_{1},\nu}.$

{\bf Lemma 7.3.} {\it Let the following conditions hold:

1) metric tensor of $R^{3}$ satisfies the conditions: $\exists
M_{0}=const>0,$ such that
$\|\tilde{a}_{\alpha\beta}\|_{m,\nu}<M_{0},$ $\|\partial
\tilde{a}_{\alpha\beta}\|_{m,\nu}<M_{0},$ $\|\partial^{2}
\tilde{a}_{\alpha\beta}\|_{m,\nu}<M_{0}.$

2) $\exists t_{0}>0,$ such that $c(t), c_{,i}(t), a^{k}(t),
\partial_{i}a^{k}(t)$ are continuous by $t, \forall
t\in[0,t_{0}],$ $c(0)\equiv 0, c_{,i}(0)\equiv 0,
 a^{k}(0)\equiv 0, \partial_{i}a^{k}(0)\equiv 0.$

3) $\exists t_{0}>0,$ such that $z^{\alpha}(t)\in C^{m-2,\nu} ,
z^{\alpha}_{,i}(t)\in C^{m-3,\nu},$ $\forall t\in[0,t_{0}].$

Then $\forall \varepsilon>0 \exists t_{0}>0$ such that

1) $ \|W_{1}\|^{(t)}_{m-3,\nu} \leq \varepsilon, \forall t\in
[0,t_{0}].$

2) $ \|W_{2}\|^{(t)}_{m-3,\nu} \leq \varepsilon, \forall t\in
[0,t_{0}].$

3) $ \|\Psi_{2}\|^{(t)}_{m-3,\nu} \leq \varepsilon, \forall t\in
[0,t_{0}].$ }

Proof of lemma follows from the form of functions $W_{1}, W_{2},
\Psi_{2},$
 properties of space $C^{m,\nu}$ and previous lemmas.

\section*{
\S 8. Properties of functions $\dot{W}_{1},\dot{W}_{2},
\dot{\Psi}_{2}.$}

\section*{
\S 8.1. Formula of function $\dot{W}_{1}.$}

We have:
$$
\dot{W}_{1}= (\dot{\tilde{a}}_{\alpha\beta}(t)- \frac{\partial
\tilde{a}_{\alpha\beta}}{\partial y^{\sigma}}(0) \dot{z}^{\sigma})
g^{ij}y^{\alpha},_{i}y^{\beta},_{j}+$$
$$
 2\dot{\tilde{a}}_{\alpha\beta}(t)g^{ij}y^{\alpha},_{i}z^{\beta},_{j}+
2(\tilde{a}_{\alpha\beta}(t)-\tilde{a}_{\alpha\beta}(0))
g^{ij}y^{\alpha},_{i}\dot{z}^{\beta},_{j}+
\dot{\tilde{a}}_{\alpha\beta}(t)g^{ij}z^{\alpha},_{i}z^{\beta},_{j}+$$
$$\tilde{a}_{\alpha\beta}(t)g^{ij}\dot{z}^{\alpha},_{i}z^{\beta},_{j}+
\tilde{a}_{\alpha\beta}(t)g^{ij}z^{\alpha},_{i}\dot{z}^{\beta},_{j}.
\eqno(8.1)$$

Consider the formula: $$\dot{\tilde{a}}_{\alpha\beta}(t)=
\frac{\partial \tilde{a}_{\alpha\beta}(t)}{\partial
y^{\sigma}}\dot{z}^{\sigma}. \eqno(8.2)$$

Using (8.2) we obtain:
$$\dot{W}_{1}= (
\frac{\partial \tilde{a}_{\alpha\beta}(t)}{\partial
y^{\sigma}}\dot{z}^{\sigma}- \frac{\partial
\tilde{a}_{\alpha\beta}}{\partial y^{\sigma}}(0) \dot{z}^{\sigma})
g^{ij}y^{\alpha},_{i}y^{\beta},_{j}+$$
$$
 2\frac{\partial \tilde{a}_{\alpha\beta}(t)}{\partial y^{\sigma}}\dot{z}^{\sigma}
 g^{ij}y^{\alpha},_{i}z^{\beta},_{j}+
2(\tilde{a}_{\alpha\beta}(t)-\tilde{a}_{\alpha\beta}(0))
g^{ij}y^{\alpha},_{i}\dot{z}^{\beta},_{j}+ \frac{\partial
\tilde{a}_{\alpha\beta}(t)}{\partial y^{\sigma}}\dot{z}^{\sigma}
g^{ij}z^{\alpha},_{i}z^{\beta},_{j}+
$$
$$
\tilde{a}_{\alpha\beta}(t)g^{ij}\dot{z}^{\alpha},_{i}z^{\beta},_{j}+
\tilde{a}_{\alpha\beta}(t)g^{ij}z^{\alpha},_{i}\dot{z}^{\beta},_{j}.
\eqno(8.3)$$

\section*{
\S 8.2. Formula of function $\dot{W}_{2}.$}

$$\dot{W}_{2}=\dot{\Delta}(g_{11}) \Delta(g_{22})+
\Delta(g_{11})
\dot{\Delta}(g_{22})-2\Delta(g_{12})\dot{\Delta}(g_{12}).
\eqno(8.4) $$

$$\dot{\Delta}(g_{ij})=
\tilde{a}_{\alpha\beta}(0)y^{\alpha},_{i}\nabla^{*}_{j}\dot{z}^{\beta}+
\tilde{a}_{\alpha\beta}(0)y^{\beta},_{j}\nabla^{*}_{i}\dot{z}^{\alpha}+
(\dot{\tilde{a}}_{\alpha\beta}(t)- \frac{\partial
\tilde{a}_{\alpha\beta}}{\partial y^{\sigma}}(0) \dot{z}^{\sigma})
y^{\alpha},_{i}y^{\beta},_{j}+$$
$$ \dot{\tilde{a}}_{\alpha\beta}(t)
 (y^{\alpha},_{i}z^{\beta},_{j}+y^{\beta},_{j}z^{\alpha},_{i})+
(\tilde{a}_{\alpha\beta}(t)-\tilde{a}_{\alpha\beta}(0))
 (y^{\alpha},_{i}\dot{z}^{\beta},_{j}+y^{\beta},_{j}\dot{z}^{\alpha},_{i})+
$$
$$\dot{\tilde{a}}_{\alpha\beta}(t)z^{\alpha},_{i}z^{\beta},_{j}+
\tilde{a}_{\alpha\beta}(t)\dot{z}^{\alpha},_{i}z^{\beta},_{j}+
\tilde{a}_{\alpha\beta}(t)z^{\alpha},_{i}\dot{z}^{\beta},_{j}.
\eqno(8.5)$$

Then we have:
$$\dot{\Delta}(g_{ij})=
\tilde{a}_{\alpha\beta}(0)y^{\alpha},_{i}\nabla^{*}_{j}\dot{z}^{\beta}+
\tilde{a}_{\alpha\beta}(0)y^{\beta},_{j}\nabla^{*}_{i}\dot{z}^{\alpha}+
(\frac{\partial \tilde{a}_{\alpha\beta}(t)}{\partial
y^{\sigma}}\dot{z}^{\sigma}- \frac{\partial
\tilde{a}_{\alpha\beta}}{\partial y^{\sigma}}(0) \dot{z}^{\sigma})
y^{\alpha},_{i}y^{\beta},_{j}+$$
$$ \frac{\partial \tilde{a}_{\alpha\beta}(t)}{\partial y^{\sigma}}\dot{z}^{\sigma}
 (y^{\alpha},_{i}z^{\beta},_{j}+y^{\beta},_{j}z^{\alpha},_{i})+
(\tilde{a}_{\alpha\beta}(t)-\tilde{a}_{\alpha\beta}(0))
 (y^{\alpha},_{i}\dot{z}^{\beta},_{j}+y^{\beta},_{j}\dot{z}^{\alpha},_{i})+
$$
$$\frac{\partial \tilde{a}_{\alpha\beta}(t)}{\partial y^{\sigma}}\dot{z}^{\sigma}
z^{\alpha},_{i}z^{\beta},_{j}+
\tilde{a}_{\alpha\beta}(t)\dot{z}^{\alpha},_{i}z^{\beta},_{j}+
\tilde{a}_{\alpha\beta}(t)z^{\alpha},_{i}\dot{z}^{\beta},_{j}.
\eqno(8.6)$$

\section*{
\S 8.3. The inequalities for norms of functions
$\dot{W}_{1},\dot{W}_{2}, \dot{\Psi}_{2}.$}

Denote:\\
$\|\nabla^{*}z \|^{(t)}_{m_{1},\nu}= \max_{\alpha,
i}\|\nabla^{*}_{i}z^{\alpha} \|^{(t)}_{m_{1},\nu}= \max_{\alpha,
i}\max_{\tau \in [0;t]}
\|\nabla^{*}_{i}z^{\alpha}(\tau)\|_{m_{1},\nu}.$\\
$\|\nabla^{*}\dot{z}\|^{(t)}_{m_{1},\nu}= \max_{\alpha,
i}\|\nabla^{*}_{i}\dot{z}^{\alpha} \|^{(t)}_{m_{1},\nu}=
\max_{\alpha, i}\max_{\tau \in [0;t]}
\|\nabla^{*}_{i}\dot{z}^{\alpha}(\tau)\|_{m_{1},\nu}.$

{\bf Lemma 8.3.1.} {\it The following estimations hold:

1) $ \|\dot{W}_{1}\|^{(t)}_{m_{1},\nu} \leq M_{20}
(\|\dot{z}\|^{(t)}_{m_{1},\nu}\|z \|^{(t)}_{m_{1},\nu}+
\|\dot{z}\|^{(t)}_{m_{1},\nu}\|\partial z \|^{(t)}_{m_{1},\nu}+
\|z\|^{(t)}_{m_{1},\nu}\|\partial \dot{z} \|^{(t)}_{m_{1},\nu}+\\
\|\dot{z}\|^{(t)}_{m_{1},\nu}(\|\partial z
\|^{(t)}_{m_{1},\nu})^{2}+ \|\partial
\dot{z}\|^{(t)}_{m_{1},\nu}\|\partial z \|^{(t)}_{m_{1},\nu}).$

2) $\|\nabla^{*}z \|^{(t)}_{m_{1},\nu}\leq M_{21} (\|\partial
z\|^{(t)}_{m_{1},\nu}+\|z \|^{(t)}_{m_{1},\nu}).$

3) $\|\nabla^{*}\dot{z} \|^{(t)}_{m_{1},\nu}\leq M_{22}
(\|\partial \dot{z}\|^{(t)}_{m_{1},\nu}+\|\dot{z}
\|^{(t)}_{m_{1},\nu}).$

4) $ \|\dot{\Delta}(g_{ij})\|^{(t)}_{m_{1},\nu} \leq M_{23}
(\|\nabla^{*}\dot{z} \|^{(t)}_{m_{1},\nu}+
\|\dot{z}\|^{(t)}_{m_{1},\nu}\|z \|^{(t)}_{m_{1},\nu}+
\|\dot{z}\|^{(t)}_{m_{1},\nu}\|\partial z\|^{(t)}_{m_{1},\nu}+
\|z\|^{(t)}_{m_{1},\nu}\|\partial \dot{z} \|^{(t)}_{m_{1},\nu}+
\|\dot{z}\|^{(t)}_{m_{1},\nu}(\|\partial z
\|^{(t)}_{m_{1},\nu})^{2}+ \|\partial
\dot{z}\|^{(t)}_{m_{1},\nu}\|\partial z \|^{(t)}_{m_{1},\nu}).$

5) $ \|\dot{W}_{2}\|^{(t)}_{m_{1},\nu} \leq M_{24}
(\max_{i,j}\|\dot{\Delta}(g_{ij})\|^{(t)}_{m_{1},\nu})
(\max_{i,j}\|\Delta(g_{ij})\|^{(t)}_{m_{1},\nu}).$ }

Proof of lemma follows from the forms of functions $W_{1}, W_{2},
\dot{W}_{1}, \dot{W}_{2},$
 properties of space $C^{m_{1},\nu}$ and previous lemmas.

{\bf Lemma 8.3.2.} {\it Let the conditions of lemma 7.3. hold:

1) metric tensor of $R^{3}$ satisfies the conditions: $\exists
M_{0}=const>0,$ such that
$\|\tilde{a}_{\alpha\beta}\|_{m,\nu}<M_{0},$ $\|\partial
\tilde{a}_{\alpha\beta}\|_{m,\nu}<M_{0},$ $\|\partial^{2}
\tilde{a}_{\alpha\beta}\|_{m,\nu}<M_{0}.$

2) $\exists t_{0}>0,$ such that $c(t), c_{,i}(t), a^{k}(t),
\partial_{i}a^{k}(t)$ are continuous by $t, \forall
t\in[0,t_{0}],$ $c(0)\equiv 0, c_{,i}(0)\equiv 0,
 a^{k}(0)\equiv 0, \partial_{i}a^{k}(0)\equiv 0.$

3) $\exists t_{0}>0,$ such that $z^{\alpha}(t)\in C^{m-2,\nu} ,
z^{\alpha}_{,i}(t)\in C^{m-3,\nu},$ $\forall t\in[0,t_{0}].$

Then $\forall \varepsilon>0 \exists t_{0}>0,$ such that

1) $ \|\dot{W}_{1}\|^{(t)}_{m-3,\nu} \leq \varepsilon, \forall
t\in [0,t_{0}].$

2) $ \|\dot{W}_{2}\|^{(t)}_{m-3,\nu} \leq \varepsilon, \forall
t\in [0,t_{0}].$

3) $ \|\dot{\Psi}_{2}\|^{(t)}_{m-3,\nu} \leq \varepsilon, \forall
t\in [0,t_{0}].$ }

Proof of lemma follows from the form of functions $W_{1}, W_{2},
\dot{W}_{1}, \dot{W}_{2},$
 properties of space $C^{m_{1},\nu}$ and previous lemmas.

\section*{\S 9. Decidability of boundary-value problem $A.$}

We have the following equation system of elliptic type:
$$ \partial_{2}\dot{a}^{1}-\partial_{1}\dot{a}^{2}+p_{k}\dot{a}^{k}=
\dot{\Psi_{1}}, $$
$$ \partial_{1}\dot{a}^{1}+\partial_{2}\dot{a}^{2}+q^{(b)}_{k}\dot{a}^{k}=
\dot{\Psi}^{(b)}_{2}, \eqno(9.1) $$ where we use (6.2.53).
$\dot{\Psi}^{(b)}_{2}=q^{(b)}_{0}\dot{c}-P_{0}.$
 Note that $q^{(b)}_{k}$ do not depend on $t.$

Without loss of generality we denote $x^{1}$ as $x^{2}$ and
$x^{2}$ as $x^{1}.$

We write (9.1) as:
$$ \partial_{1}\dot{a}^{1}-\partial_{2}\dot{a}^{2}+p_{k}\dot{a}^{k}=
\dot{\Psi_{1}}, $$
$$ \partial_{2}\dot{a}^{1}+\partial_{1}\dot{a}^{2}+q^{(b)}_{k}\dot{a}^{k}=
\dot{\Psi}^{(b)}_{2}, \eqno(9.2) $$

Define: $w=a^{1}+ia^{2}, z=x^{1}+i x^{2}.$

Therefore we have boundary-value problem for generalized analytic
functions.

$$ \partial_{\bar{z}}\dot{w}+A\dot{w}+B\bar{\dot{w}}=\dot{\Psi_{0}},
\qquad Re\{\overline{\lambda}\dot{w}\}=\dot{\varphi} \quad on
\quad \partial D, \eqno(9.3) $$ where
$$\partial_{\bar{z}}\dot{w}=\frac{1}{2}(\dot{w}_{x}+i
\dot{w}_{y}), A=\frac{1}{4}(p_{1}+q^{(b)}_{2}+i q^{(b)}_{1}-i
p_{2}),$$
$$B=\frac{1}{4}(p_{1}-q^{(b)}_{2}+i q^{(b)}_{1}+i p_{2}),
\dot{\Psi}_{0}=\frac{1}{2}(\dot{\Psi}_{1}+i \dot{\Psi}^{(b)}_{2}).
\eqno(9.4)$$

We change the form of obtained boundary-value problem (9.4).
Consider the following:

$$\dot{\Psi}^{(b)}_{2}=q^{(b)}_{0}\dot{c}-P_{0}=$$
$$ q^{(b)}_{0}\Biggl(\int\limits_{(x^{1}_{(0)},x^{2}_{(0)})}^{(x^{1},x^{2})}
\Biggl(-V\dot{a}^{1}\Biggr) d\tilde{x}^{1}+
\Biggl(-V\dot{a}^{2}\Biggr)d\tilde{x}^{2}\Biggr) +q^{(b)}_{0}
P(\dot{a}^{1},\dot{a}^{2})-P_{0}. \eqno(9.5) $$

Denote:
$$\dot{\Psi}_{3}= q^{(b)}_{0} P(\dot{a}^{1},\dot{a}^{2})
 -P_{0}. \eqno(9.6) $$

We define: $\dot{\Psi}=\frac{1}{2}(\dot{\Psi}_{1}+i
\dot{\Psi}_{3}).$

By (9.5), (9.6), the boundary-value problem (9.4) takes the form:
$$\partial_{\bar{z}}\dot{w}+A\dot{w}+B\bar{\dot{w}}+
i
\frac{q^{(b)}_{0}}{2}\int\limits_{(x^{1}_{(0)},x^{2}_{(0)})}^{(x^{1},x^{2})}
\Biggl(V\dot{a}^{1}\Biggr) d\tilde{x}^{1}+
\Biggl(V\dot{a}^{2}\Biggr)d\tilde{x}^{2}=\dot{\Psi}, \eqno(9.7) $$
$Re\{\overline{\lambda}\dot{w}\}=\dot{\varphi}$ on $\partial D.$

Therefore we have:
$$\partial_{\bar{z}}\dot{w}+A\dot{w}+B\bar{\dot{w}}+$$
$$\int\limits_{(x^{1}_{(0)},x^{2}_{(0)})}^{(x^{1},x^{2})}
\Biggl(i
\frac{q^{(b)}_{0}(x^{1},x^{2})}{2}V(\tilde{x}^{1},\tilde{x}^{2})
\dot{a}^{1}\Biggr) d\tilde{x}^{1}+ \Biggl(i
\frac{q^{(b)}_{0}(x^{1},x^{2})}{2}V(\tilde{x}^{1},\tilde{x}^{2})
\dot{a}^{2}\Biggr)d\tilde{x}^{2}=\dot{\Psi}, \eqno(9.8) $$
$Re\{\overline{\lambda}\dot{w}\}=\dot{\varphi}$ on $\partial D.$

Using the formulas:
$\dot{a}^{1}=\frac{1}{2}(\dot{w}+\bar{\dot{w}}),$
$\dot{a}^{2}=\frac{i}{2}(\bar{\dot{w}}-\dot{w}),$\\
and denoting: $E_{0}=i
\frac{q^{(b)}_{0}(x^{1},x^{2})}{2}V(\tilde{x}^{1},\tilde{x}^{2})$
we obtain the following form of desired boundary-value problem:
$$\partial_{\bar{z}}\dot{w}+A\dot{w}+B\bar{\dot{w}}+
\int\limits_{(x^{1}_{(0)},x^{2}_{(0)})}^{(x^{1},x^{2})}
\Biggl(\frac{E_{0}}{2}(\dot{w}+\bar{\dot{w}})  \Biggr)
d\tilde{x}^{1}+ \Biggl(\frac{i E_{0}}{2}(\bar{\dot{w}}-\dot{w})
\Biggr)d\tilde{x}^{2}=\dot{\Psi}, \eqno(9.9) $$
$Re\{\overline{\lambda}\dot{w}\}=\dot{\varphi}$ on $\partial D.$

We denote:
$$E(\dot{w})=\int\limits_{(x^{1}_{(0)},x^{2}_{(0)})}^{(x^{1},x^{2})}
\Biggl(\frac{E_{0}}{2}(\dot{w}+\bar{\dot{w}})  \Biggr)
d\tilde{x}^{1}+ \Biggl(\frac{i E_{0}}{2}(\bar{\dot{w}}-\dot{w})
\Biggr)d\tilde{x}^{2}. \eqno(9.10)$$

Then we finally have the form of desired boundary-value problem:
$$\partial_{\bar{z}}\dot{w}+A\dot{w}+B\bar{\dot{w}}+E(\dot{w})=\dot{\Psi},
\qquad Re\{\overline{\lambda}\dot{w}\}=\dot{\varphi} \quad on
\quad \partial D. \eqno(9.11) $$

Let, along the $\partial F$, be given vector field tangent to $F.$
We denote it by the following formula:
$$v^{\alpha}=l^{i}y^{\alpha}_{,i}.  \eqno(9.12)$$

We consider the boundary-value condition:
$$\tilde{a}_{\alpha\beta}z^{\alpha}v^{\beta}=\tilde{\gamma}(s,t) ,
s\in \partial D. \eqno(9.13)$$

Define:
$\tilde{\lambda}_{k}=\tilde{a}_{\alpha\beta}y^{\alpha}_{,k}v^{\beta},
k=1,2.$ \\ Then boundary condition takes the form:
$Re\{(\dot{a^{1}}+i \dot{a^{2}}) (\tilde{\lambda}_{1}-i
\tilde{\lambda}_{2})\}=\dot{\tilde{\gamma}}$ on $\partial F.$

Denote: $\lambda_{k}=\frac{ \tilde{\lambda}_{k}}
{(\tilde{\lambda}_{1})^{2}+(\tilde{\lambda}_{2})^{2}}, k=1,2.$
$\dot{\varphi}=\frac{ \dot{\tilde{\gamma}} }
{(\tilde{\lambda}_{1})^{2}+(\tilde{\lambda}_{2})^{2}}.$

Then boundary-value condition takes the form:
$Re\{\overline{\lambda} \dot{w}\}=\dot{\varphi}$ on $\partial F,$
where $|\lambda|=1.$

We analyze the decidability of the following boundary-value problem (A):\\
$$\partial_{\bar{z}}\dot{w}+A\dot{w}+B\bar{\dot{w}}+E(\dot{w})=\dot{\Psi},
\qquad Re\{\overline{\lambda}\dot{w}\}=\dot{\varphi} \quad on
\quad \partial D, \eqno(9.14) $$
$\lambda=\lambda_{1}+i\lambda_{2},$ $|\lambda|\equiv 1,$ $\lambda,
\dot{\varphi}\in C^{m-2,\nu}(\partial D).$

We use the fact that $\dot{\Psi}=\dot{\Psi}(\dot{w},z,t),
E(\dot{w})=E(\dot{w},z,t),$
$\dot{w}=\dot{w}(t),$ \\
$\dot{\varphi}=\dot{\varphi}(s,t), s\in \partial D,$
$\lambda=\lambda(s), s\in \partial D.$

Let $n$ be index of obtained boundary-value problem
$$n=\frac{1}{2\pi}\Delta_{\partial D} \arg \lambda(s). \eqno(9.15) $$

{\bf Theorem 9.1.} {\it Let $t$ be fixed.\\
Let $A(z), B(z), \dot{\Psi}(z) \in C^{m-3,\nu}(\bar{D}),$
$\lambda(s), \dot{\varphi}\in C^{m-2,\nu}(\partial D),$
$|\lambda(s)|\equiv 1.$

Let $\dot{\Psi}(0,z)=0,$
$\|\dot{\Psi}(\dot{w_{1}},z)-\dot{\Psi}(\dot{w_{2}},z)
\|_{m-2,\nu}\leq
\mu(\rho)\|\dot{w_{1}}-\dot{w_{2}}\|_{m-1,\nu},$ \\
for $\|\dot{w}_{1}\|_{m-2,\nu}\leq \rho,
\|\dot{w}_{2}\|_{m-2,\nu}\leq \rho ,$ $\lim_{\rho \rightarrow 0}
\mu(\rho) =0.$

Then, assuming that $t$ is fixed, the following holds:

1) if $n\geq 0$ then there exist $\rho$ and $\varepsilon(\rho)>0$
such that for $\|\dot{\varphi}\|_{m-2,\nu}\leq \varepsilon$ the
boundary-value problem has $(2n+1)-$parametric solution of class
$C^{m-2,\nu}(\bar{D})$ for any admissible $\dot{\varphi}$.

2) if $n < 0$ then there exist $\rho>0$ and $\varepsilon(\rho)>0$
such that for $\|\dot{\varphi}\|_{m-2,\nu}\leq \varepsilon(\rho)$
the boundary-value problem has nor more than one solution of class
$C^{m-2,\nu}(\bar{D})$ for any admissible $\dot{\varphi}$. For
$\dot{\varphi}\equiv 0$ boundary-value problem with condition:
$\|\dot{w}\|_{m-2,\nu}\leq \rho$ has only zero solution. }

{\bf Proof.} Consider the following boundary-value problem $(A_{0})$:\\
$$\partial_{\bar{z}}\dot{w}+A\dot{w}+B\bar{\dot{w}}=\dot{\Psi}, \qquad
Re\{\overline{\lambda}\dot{w}\}=\dot{\varphi} \quad on \quad
\partial D, \eqno(9.16)$$
 $|\lambda|\equiv 1,$
$\dot{\varphi}\in C^{m-2,\nu}(\partial D)$, $\dot{\Psi}\in
C^{m-3,\nu}(\bar{D}).$

Consider the operator:\\
$$
I(\dot{\Psi},z)=-\frac{1}{\pi}\int\int\limits_{D}(\Omega_{1}(z,\zeta)\dot{\Psi}(\zeta)+
\Omega_{2}(z,\zeta)\overline{\dot{\Psi}(\zeta)})d\xi d\eta,
\zeta=\xi+i\eta , \eqno(9.17)$$ where $\Omega_{1}, \Omega_{2}$ are
principal kernels of the equation
$\partial_{\bar{z}}\dot{w}+A(z)\dot{w}+B(z)\bar{\dot{w}}=0.$

It is well known [17,18] that operator $I(\dot{\Psi},z)$ takes the form:\\
$$
I(\dot{\Psi},z)=T(\dot{\Psi})-\frac{1}{\pi}\int\int\limits_{D}(K_{1}(z,\zeta)\dot{\Psi}(\zeta)+
\Omega_{2}(z,\zeta)\overline{\dot{\Psi}(\zeta)})d\xi d\eta,
\zeta=\xi+i\eta , \eqno(9.18)$$
$$
T(\dot{\Psi})=-\frac{1}{\pi}\int\int\limits_{D}
\frac{\dot{\Psi}(\zeta)}{\zeta-z}d\xi d\eta, \eqno(9.19) $$ where
operator $T(\dot{\Psi})$ is completely continuous [17,18].

Consider the operator:\\
$$A(\dot{\Psi},z)=I(\dot{\Psi},z)+
\int\limits_{\partial
D}Re\{\overline{\lambda(s)}I(\dot{\Psi},s)\}M_{0}(z,s)ds,
\eqno(9.20)
$$
where $M_{0}(z,s)$ is kernel of boundary-value problem (do not
depend on $\dot{\varphi}$)
$$\partial_{\bar{z}}\dot{w}+A(z)\dot{w}+B(z)\bar{\dot{w}}=0,
Re\{\overline{\lambda(s)}w(s)\}=\dot{\varphi}, s \in \partial D.
\eqno(9.21) $$

Consider the operator
$$A_{2}(\dot{w})=A_{1}(\dot{w})=A(\dot{\Psi}(\dot{w},z)). \eqno(9.22)$$

According to the results from [20], theorem 9.1. is valid for
problem $(A_{0})$.
For the case $n\geq 0$ problem $(A_{0})$ is solved as:\\
$$
\dot{w}=A_{2}(\dot{w})+\int\limits_{\partial D}
\dot{\varphi}(s)M_{0}(z,s)ds+ \sum_{i=1}^{2n+1}c_{i}\dot{w}_{i}.
\eqno(9.23)$$

Therefore for the case $n\geq 0$ problem $(A)$ is solved as:\\
$$\dot{w}=A_{2}(\dot{w})+
\int\limits_{\partial D} \dot{\varphi}(s)M_{0}(z,s)ds+
\sum_{i=1}^{2n+1}c_{i}\dot{w}_{i}+A_{2}(E(\dot{w})). \eqno(9.24)$$
Then for this equation we use theory of Fredholm operator of index
zero and theory of Volterra operator equation. Therefore we can
solve (9.24) by the method of successive approximations.

For the case $n < 0$ we solve the problem $(A_{0})$ as equation
system
consisting of $-2n$ equations:\\
$$
\dot{w}=A_{2}(\dot{w})+\int\limits_{\partial D}
\dot{\varphi}(s)M_{0}(z,s)ds, \eqno(9.25) $$
$$\int\limits_{\partial D} (\dot{\varphi}(s)+
Re\{\overline{\lambda(s)}I(\dot{\Psi},s)\})\dot{w}'_{j}(s)\lambda(s)ds=0,
\quad j=\overline{1,-2n-1},$$ where $\dot{w}'_{j}$ are complete
system of solutions of the following problem:
$$\partial_{\bar{z}}\dot{w}'-A(z)\dot{w}'-\bar{B}(z)
\overline{\dot{w}}'=0, \quad Re\{\lambda(z)\frac{d z(s)}{d s}
\dot{w}'(z)\}=0 \quad on \quad \partial D
$$

Then for the case $n < 0$ we solve the problem $(A)$ as equation
system
consisting of $-2n$ equations:\\
$$
\dot{w}=A_{2}(\dot{w})+\int\limits_{\partial D}
\dot{\varphi}(s)M_{0}(z,s)ds+ A_{2}(E(\dot{w})), \eqno(9.26) $$
$$\int\limits_{\partial D} (\dot{\varphi}(s)+
Re\{\overline{\lambda(s)}I(\dot{\Psi},s)\})\dot{w}'_{j}(s)\lambda(s)ds=0,
\quad j=\overline{1,-2n-1},$$ where $\dot{w}'_{j}$ are complete
system of solutions of the following problem:
$$\partial_{\bar{z}}\dot{w}'-A(z)\dot{w}'-\bar{B}(z)
\overline{\dot{w}}'=0, \quad Re\{\lambda(z)\frac{d z(s)}{d s}
\dot{w}'(z)\}=0 \quad on \quad \partial D
$$
Then for this equation system we use theory of Fredholm operator
of index zero and theory of Volterra operator equation.

By modifying standard method from [20], using the method of
successive approximations and principle of contractive mapping, we
obtain the proof theorem 9.1. for boundary-value problem $(A).$

{\bf Theorem 9.2.} {\it Let $F\in C^{m,\nu}, \nu \in (0;1) , m\ge
4,$ $\partial F\in C^{m+1,\nu}.$

Then the following holds:

1) if $n\geq 0$ then there exists $t_{0}>0$ and exists
$\varepsilon(t_{0})>0$ such that for
$\|\dot{\varphi}\|_{m-2,\nu}\leq \varepsilon$ boundary-value
problem (A) for all $t\in [0, t_{0})$ has $(2n+1)-$parametric
solution of class $C^{m-2,\nu}(\bar{D})$ continuous by $t\in
[0,t_{0})$ for any admissible $\dot{\varphi}.$

2) if $n < 0$ then exists $t_{0}>0$ and exists
$\varepsilon(t_{0})>0$ such that for
$\|\dot{\varphi}\|_{m-2,\nu}\leq \varepsilon(t_{0})$
boundary-value problem (A) for all $t\in [0, t_{0})$ has nor more
than one solution of class $C^{m,\nu}(\bar{D})$ continuous by
$t\in [0,t_{0})$ for any admissible $\dot{\varphi}$. For
$\dot{\varphi}\equiv 0$ the boundary-value problem has only zero
solution. }

Proof follows from theorem 9.1., form of function $\dot{\Psi}$ and
the fact that for all sufficiently small $t$ the conditions of
theorem 9.1 hold.

\section*{\S 10. Proof of theorem 1.}

Proof of theorem 1 follows from theorem 9.2., formulas of
$MG-$deformation and formulas of finding function $\dot{c}$ on
functions $\dot{a}^{j}.$ Using the condition of theorem 1: at the
point $(x^{1}_{(0)},x^{2}_{(0)})$ of the domain $D,$ the following
condition holds: $\forall t : a^{i}(t)\equiv 0, c(t)\equiv 0.$
Therefore in case 1) $n > 0$ boundary-value problem (A) has
$(2n-1)-$parametric solution. Using similar reasonings we prove
theorem 1 for the cases 2) and 3).

The theorem 1 is proved.

\section*{\bf References.}

\begin{enumerate}

 \item A.I. Bodrenko.
On continuous almost ARG-deformations of hypersurfaces in
Euclidean space [in Russian]. Dep. in VINITI 27.10.92., N3084-T92,
UDK 513.81, 14 pp.

 \item A.I. Bodrenko. Some properties continuous ARG-deformations
 [in Russian]. Theses of international science conference
"Lobachevskii and modern geometry", Kazan, Kazan university
publishing house, 1992 ., pp.15-16.

 \item A.I. Bodrenko. On continuous ARG-deformations [in Russian].
Theses of reports on republican science and methodical conference,
dedicated to the 200-th anniversary of N.I.Lobachevskii, Odessa,
Odessa university publishing house, 1992 ., Part 1, pp.56-57.

 \item A.I. Bodrenko. On extension of infinitesimal almost
 ARG-deformations closed \\
 hypersurfaces into analytic
 deformations in Euclidean spaces [in Russian].
 Dep. in VINITI 15.03.93., N2419-T93 UDK 513.81, 30 pp.

 \item A.I. Bodrenko. On extension of infinitesimal almost
 ARG-deformations of hypersurface with boundary into analytic deformations
[in Russian]. Collection works of young scholars of VolSU,
Volgograd, Volgograd State University publishing house, 1993,\\
pp.79-80.

 \item A.I. Bodrenko. Some properties of continuous almost
AR-deformations of hypersurfaces with prescribed change of
Grassmannian image  [in Russian].
 Collection of science works of young scholars, Taganrog,
 Taganrog State Pedagogical Institute publishing house , 1994, pp. 113-120.

 \item A.I. Bodrenko. On continuous almost AR-deformations
with prescribed change of Grassmannian image [in Russian].
All-Russian school-colloquium on stochastic \\ methods of geometry
and analysis. Abrau-Durso. Publisher Moscow: "TVP". Theses of
reports, 1994, pp. 15-16.

 \item A.I. Bodrenko. Extension of infinitesimal almost ARG-deformations
of \\ hypersurfaces into analytic deformations [in Russian].
All-Russian school-colloquium on stochastic methods. Yoshkar-Ola.
Publisher Moscow: "TVP". Theses of reports, 1995, pp. 24-25.

 \item A.I. Bodrenko. Areal-recurrent deformations of hypersurfaces
 preserving Grassmannian image [in Russian].
Dissertation of candidate of physical-mathematical sciences. \\
Novosibirsk,1995, pp. 85.

 \item A.I. Bodrenko. Areal-recurrent deformations of hypersurfaces
preserving Grassmannian image [in Russian]. Author's summary of
dissertation of candidate of physical- \\ mathematical sciences.
Novosibirsk, 1995, pp. 1-14.

 \item A.I. Bodrenko. Some properties of ARG-deformations [in Russian].
Izvestiay Vuzov. Ser. Mathematics, 1996, N2, pp.16-19.

 \item A.I. Bodrenko. Continuous almost ARG-deformations of surfaces
 with boundary [in Russian].
Modern geometry and theory of physical fields. \\ International
geometry seminar of N.I.Lobachevskii Theses of reports, Kazan,
\\
Publisher Kazan university, 1997, pp.20-21.

 \item A.I. Bodrenko. Continuous almost AR-deformations of surfaces
with prescribed change of Grassmannian image [in Russian]. Red.
"Sib. mat. zhurnal.", Sib. otd. RAN , Novosibirsk, Dep. in VINITI
13.04.98., N1075-T98 UDK 513.81, 13 pp.

 \item A.I. Bodrenko. Almost ARG-deformations of the second order of
 surfaces in \\ Riemannian space [in Russian].
 Surveys in Applied and Industrial Mathematics. \\ 1998,
 Vol. 5, Issue 2, p.202. Publisher Moscow: "TVP".

 \item A.I. Bodrenko. Almost $AR$-deformations of a surfaces with prescribed
change of \\ Grassmannian image with exterior connections [in
 Russian]. Red. zhurn. "Izvestya vuzov. Mathematics.", Kazan, Dep.
in VINITI  03.08.98, N2471 - B 98. P. 1-9.

 \item A.I. Bodrenko. Properties of generalized G-deformations with
areal condition of normal type in Riemannian space [in Russian].
 Surveys in Applied and Industrial Mathematics.
Vol. 7. Issue 2. (VII All-Russian school-colloquium on stochastic
methods. Theses of reports.) P. 478. Moscow: TVP, 2000.

 \item I.N. Vekua. Generalized Analytic Functions. Pergamon.
 New York. 1962.

 \item I.N. Vekua. Generalized Analytic Functions [in Russian].
 Moscow. Nauka. 1988.

 \item I.N. Vekua. Some questions of the theory of differential equations
 and  applications in mechanics [in Russian].
 Moscow:"Nauka". 1991 . pp. 256.

 \item A.V. Zabeglov. On decidability of one nonlinear
  boundary-value problem for AG-deformations
  of surfaces with boundary [in Russian].
  Collection of science works. \\ Transformations of surfaces,
  Riemannian spaces determined by given recurrent \\ relations. Part 1.
Taganrog. Taganrog State Pedagogical Institute publishing house.
1999. pp. 27-37.

 \item V.T. Fomenko. On solution of the generalized Minkowski problem
 for surface with boundary [in Russian].
Collection of science works. Transformations of surfaces,\\
  Riemannian spaces determined by given recurrent relations. Part 1.
Taganrog. \\ Taganrog State Pedagogical Institute publishing
house. 1999. pp. 56-65.

 \item V.T. Fomenko. On uniqueness of solution of the generalized
 Christoffel problem for surfaces with boundary [in Russian].
Collection of science works. Transformations of surfaces,
  Riemannian spaces determined by given recurrent relations. Part 1.
Taganrog. Taganrog State Pedagogical Institute publishing house.
1999. pp. 66-72.

 \item V.T. Fomenko. On rigidity of surfaces with boundary in
 Riemannian space [in Russian].
 Doklady Akad. Nauk SSSR. 1969 . Vol. 187, N 2, pp. 280-283.

 \item V.T. Fomenko. ARG-deformations of hypersurfaces in Riemannian
space [in Russian].\\  //Dep. in VINITI 16.11.90 N5805-B90

  \item S.B. Klimentov. On one method of construction the solutions
  of boundary-value \\ problems in the bending theory of surfaces of positive
  curvature [in Russian]. \\ Ukrainian geometry sbornik. pp. 56-82.

 \item M.A. Krasnoselskii. Topological methods in the theory of
 nonlinear problems

   [in Russian]. Moscow, 1965.

 \item L.P. Eisenhart. Riemannian geometry [in Russian]. \\ Izd. in. lit.,
 Moscow 1948.
 (Eisenhart Luther Pfahler. Riemannian geometry. 1926.)

 \item J.A. Schouten, D.J. Struik. Introduction into new methods
 of differential geometry [in Russian]. Volume 2. Moscow. GIIL. 1948 .
 (von J.A. Schouten und D.J. Struik.
 Einf$\ddot{u}$hrung in die neueren methoden der \\
 differentialgeometrie.
 Zweite vollst$\ddot{a}$ndig umgearbeitete
 Auflage. Zweiter band. 1938. )

 \item I. Kh. Sabitov. //VINITI. Results of science and technics.
 Modern problems of \\mathematics [in Russian].
 Fundamental directions. Vol.48, pp.196-271.

\end{enumerate}

\end{document}